\documentclass{svproc}

\usepackage[utf8]{inputenc}
\usepackage[T1]{fontenc}
\usepackage{lmodern}
\usepackage{microtype}
\usepackage{amsmath,amssymb}
\usepackage{mathdots}
\usepackage{ccicons}

\usepackage{mathtools}
\def\mystackrel#1#2{\stackrel{\scalebox{.5}{(#1)}}{#2}}

\usepackage{stmaryrd}
\def\iverson#1{\llbracket#1\rrbracket}
\def\iversonzero#1{\iverson{#1 = 0}}

\usepackage{comment}

\usepackage{url}

\usepackage{xcolor}

\def\diag{\operatorname{diag}}
\def\even{\operatorname{even}}
\def\res{\operatorname{res}}

\def\bQ{\mathbb{Q}}
\def\bZ{\mathbb{Z}}

\def\cA{\mathcal{A}}
\def\cD{\mathcal{D}}
\def\cG{\mathcal{G}}

\def\cL{\mathcal{L}}
\def\cP{\mathcal{P}}
\def\cR{\mathcal{R}}
\def\cU{\mathcal{U}}

\newcommand{\invfact}[1]{#1\text{!`}}

\usepackage[plainpages=false,pdfpagelabels,colorlinks=true,allcolors=blue!80!black,hypertexnames=false,pdftex=true,unicode]{hyperref}
\usepackage[style=numeric, maxnames=50, sorting=nyt, isbn=false, uniquename=false]{biblatex}
\begin{document}
\mainmatter
\title{On $(3,1)$-regular graphs with one more vertex than edges}
\author{Frédéric Chyzak\inst{1} \and Hui Huang\inst{2} \and Manuel Kauers\inst{3}}
\institute{%
Inria, Palaiseau, France
\and
Fuzhou University, Fuzhou, China
\and
Johannes Kepler University, Linz, Austria
}
\maketitle

\begin{abstract}
Sequence A339987 of the OEIS 
counts $(3,1)$-regular graphs having one more vertex than edges by half the number of vertices.
A~recurrence relation satisfied by this sequence
was guessed by Kauers and Koutschan in~2023.
We confirm it in three ways:
first, by a representation as the diagonal of a triple sum
and an elaborate variant of traditional creative telescoping
that makes an a~posteriori validation possible;
second, by a residue representation
and a direct calculation by reduction-based creative telescoping;
third, by a combinatorial recurrence on graph families
and a calculation by differential elimination.
Each of those three approaches leads to a formally complete proof
and involves a computer calculation in one way or another.
\end{abstract}

\section{Introduction}
\label{sec:intro}

The goal of the present work is to address a conjecture by Kauers and Koutschan~\cite{KauersKoutschan-2023-SDF},
and to prove it in the form of the following theorem.

\begin{theorem}
Sequence A339987 of the OEIS \cite{OEIS}, which
counts $(3,1)$-regular vertex-labelled graphs having $2k$~vertices and $2k-1$~edges,
satisfies the following recurrence relation,
valid for all~$k\geq0$:
\begin{equation}\label{eq:conj-rec}
\begin{aligned}
{} & 32(328k^3 + 3300k^2 + 10844k + 11589) \\
  &\qquad{}\times (k+1)(k+2)(2k+1)(2k+3)(2k+5)(2k+7)(2k+9) \, a_k \\
{}&- 8(2624k^4 + 30664k^3 + 129460k^2 + 232328k + 148119) \\
  &\qquad{}\times (k+2)(2k+3)(2k+5)(2k+7)(2k+9) \, a_{k+1} \\
{}&- 16(2952k^5 + 40852k^4 + 219308k^3 + 569267k^2 + 712135k + 341634) \\
  &\qquad{}\times (k+3)(2k+5)(2k+7)(2k+9) \, a_{k+2} \\
{}&+ 8(3936k^5 + 55672k^4 + 306380k^3 + 818282k^2 + 1057879k + 527520) \\
  &\qquad{}\times (k+4)(2k+7)(2k+9) \, a_{k+3} \\
{}&- 2(2624k^5 + 42472k^4 + 264028k^3 + 786236k^2 + 1117119k + 601452) \\
  &\qquad{}\times (k+5)(2k+9) \, a_{k+4} \\
{}&+ 3(328k^3 + 2316k^2 + 5228k + 3717) (k+4)(k+6) \, a_{k+5} = 0 .
\end{aligned}
\end{equation}
\end{theorem}
The 5-order recurrence relation~\eqref{eq:conj-rec} was announced
and formulated in a renormalized form as Conjecture~16 in~\cite{KauersKoutschan-2023-SDF},
and it can be found in the form above in the OEIS.
We will prove it formally
by three independent methods in Sections \ref{sec:triplesum}, \ref{sec:residues}, and~\ref{sec:graphs}.
First, in Section~\ref{sec:triplesum} we represent~$a_k$ by a triple sum~\eqref{eq:triple-sum-with-natural-bounds}
to which we apply the creative-telescoping algorithmic theory
\cite{Zeilberger-1991-MCT,Chyzak-2000-EZF,Koutschan-2010-FAC,BrochetSalvy-2024-RBC}.
While this outputs the correct recurrence relation~\eqref{eq:conj-rec} for the sequence~$(a_k)_k$,
the three layers of certificates that come with it
lead to too many, too singular, and too unappealing auxiliary sums to be computed,
and we were unable to complete a proof by using this output directly.
Instead, we had to resort to an expert use of a variant creative-telescoping algorithm
and to reshape the so obtained certificate
so as to make an a posteriori validation of the result possible.
This involves an unexpected use of integer linear programming to complete the formal validation.
On the other hand, the formal-residue representation that we develop in~Section~\ref{sec:residues}
lends itself to the differential counterpart of creative telescoping,
which, in the case of formal residues, is beyond any doubt a direct complete proof.
Finally, it was to be expected that Read's method of counting 3-regular graphs~\cite{Read-1970-SUE},
later simplified in the work of Wormald~\cite{Wormald-1979-ELG},
should extend to our situation.
Despite the added technicalities,
this is what we could develop into a third formal proof in Section~\ref{sec:graphs}.

The sequence~$(a_k)_{k\geq0}$ in the theorem starts with the values
\begin{multline}\label{eq:first-values}
a_0=0 , \quad a_1=1 , \quad a_2=4 , \quad a_3=90 , \quad a_4=8400 , \quad a_5=1426950 , \\
a_6=366153480 , \quad a_7=134292027870 .
\end{multline}
The rapid growth of the sequence reflects that it counts vertex-labelled graphs.
The numbers~$a_k$ were first considered by Kaygun~\cite{Kaygun-2021-ELG},
who was interested in counting the vertex-labelled graphs that realize the same degree sequences as vertex-labelled (rooted) binary trees:
they are what we call vertex-labelled $(3,1)$-regular graphs,
using a terminology~\cite{ChyzakMishna-2025-DES} that generalizes
the long-established concept of graphs all of whose vertices share the same degree~$\ell$,
usually called $\ell$-regular graphs~\cite{Read-1959-ELR,Read-1960-ELR}.
While a tree with $k+1$~leaves and therefore $k-1$~internal vertices of degree~3
is by definition connected,
general graphs with the same numbers of leaves and cubic vertices are in general not connected.
But both in the case of trees and in the case of more general graphs,
the number of edges is fixed by the degree sequence:
a graph on $2k$~vertices having $k-1$~vertices of degree~3 and $k+1$~vertices of degree~1
has $(3\times(k-1)+1\times(k+1)) / 2 = 2k-1$ edges.
Therefore, $a_k$~counts the number of vertex-labelled $(3,1)$-regular graphs having one more vertex,~$2k$, than edges,~$2k-1$.

Our proof by a triple-sum representation (Section~\ref{sec:triplesum})
and our proof by residue representation (Section~\ref{sec:residues})
are based on obtaining, respectively,
the ordinary generating function of the sequence~$(a_k)_{k\geq0}$
and a related exponential generating function.
Both are obtained
as suitable subseries extractions by scalar products in the theory of symmetric functions~\cite{Gessel-1990-SFP},
and we now introduce the relevant objects.
It will prove useful to generalize the numbers~$a_k$,
and we let $s_{m,n}$ denote the number of $(3,1)$-regular graphs on $n$~labelled vertices with $m$~unlabelled edges,
so that~$a_k = s_{2k-1,2k}$.
Furthermore, let
\begin{equation}\label{eq:Sqt-and-Hqt}
S(q, t) = \sum_{m\geq0, \ n\geq0} \frac{s_{m,n}}{n!} q^m t^n
\qquad\text{and}\qquad
H(q, t) = \sum_{m\geq0, \ n\geq0} s_{m,n} q^m t^n
\end{equation}
denote the exponential and ordinary generating functions of $(3,1)$-regular vertex-labelled graphs,
with vertices marked by~$t$ and edges marked by~$q$.
Gessel's theory~\cite{Gessel-1990-SFP},
which was later turned into algorithms~\cite{ChyzakMishnaSalvy-2005-ESP,ChyzakMishna-2025-DES},
can be used to prove the formula
\begin{equation}\label{eq:exp-gf-as-scalp}
S(q, t) = \langle \exp(f(q, p)), \exp(tg(p)) \rangle
\end{equation}
for the exponential generating function
(see the argument below).
In~\eqref{eq:exp-gf-as-scalp}, $f$ and~$g$ are two polynomials,
\begin{equation*}
f(q, p) = \tfrac12 q (p_1^2 - p_2) - \tfrac14 q^2 p_2^2 + \tfrac16 q^3 p_3^2 ,
\quad
g(p) = h_3 + h_1 = p_1 + \tfrac16 p_1^3 + \tfrac12 p_1 p_2 + \tfrac13 p_3 ,
\end{equation*}
and the scalar product~$\langle{\cdot},{\cdot}\rangle$
is the classical operation of the theory of symmetric functions,
defined by bilinearity from the formula
\begin{equation*}
\langle p_1^{s_1}\dotsm p_\ell^{s_\ell}, p_1^{r_1}\dotsm p_\ell^{r_\ell} \rangle =
  \begin{cases}
  1^{r_1} r_1!\, 2^{r_2} r_2!\, \dotsm \ell^{r_\ell} r_\ell! & \text{if $r_1 = s_1$, \dots, $r_\ell = s_\ell$,} \\
  0 & \text{otherwise.}
  \end{cases}
\end{equation*}
Proving~\eqref{eq:exp-gf-as-scalp} is a simple extension
of the theory developed in the cited references,
which deal with the case~$q=1$.
Indeed, a special case of Lemma~3.1 in~\cite{ChyzakMishna-2025-DES} is
\begin{equation}\label{eq:St}
S(t) = \langle \exp(f(p)), \exp(tg(p)) \rangle ,
\end{equation}
where $S(t) = S(1, t)$ counts graphs only by numbers of vertices
and where
\begin{equation*}
f(p) = f(1, p) = \tfrac12 p_1^2 - \tfrac14 p_2^2 + \tfrac16 p_3^2 - \tfrac12 p_2 .
\end{equation*}
The refined identity~\eqref{eq:exp-gf-as-scalp} follows immediately
after noting that $f(q, p)$~is just the homogenized version of~$f(p)$
obtained by replacing each~$p_i$ with~$q^{i/2}p_i$.

As an aside, we add that
the reader interested in a more intuitive proof will contemplate that
the (classical) reasoning of~\cite{ChyzakMishna-2025-DES} introduces
the generating function $\prod_{1\leq i < j}(1+x_ix_j)$ of all simple graphs on vertices 1, 2, \dots.
By replacing it
with the similar generating function $\prod_{1\leq i < j}(1+qx_ix_j)$,
in which each copy of~$q$ marks an edge between the vertices $i$ and~$j$,
we obtain the refined formula~\eqref{eq:exp-gf-as-scalp}.
Indeed, multiplying the product~$x_ix_j$ with a~$q$ is equivalent to multiplying each~$x_\ell$, $\ell\geq1$, with~$q^{1/2}$,
which results in multiplying~$p_\ell$,
interpreted as $\sum_{i\geq0} x_i^\ell$ in the symmetric-function theory,
with~$q^{\ell/2}$ for each~$\ell\geq1$.

The exponential generating function~$S(q, t)$ will be used in Section~\ref{sec:residues},
but in Section~\ref{sec:triplesum} we will instead consider
the ordinary generating function~$H(q, t)$ in~\eqref{eq:Sqt-and-Hqt},
with a modified scalar-product representation:
extracting the coefficient of~$t^n$ in~\eqref{eq:exp-gf-as-scalp}
and regenerating~$H(q, t)$ yields
\begin{equation}\label{eq:ord-gf-as-scalp}
H(q, t) = \left\langle \exp(f(q, p)), \frac{1}{1-tg(p)} \right\rangle .
\end{equation}
Because our sequence of interest satisfies~$a_k = s_{2k-1,2k}$,
it is the subsequence at even indices of the sequence~$(s_{n-1,n})_{n\geq0}$,
which itself is a diagonal.
More explicitly, we get
\begin{equation}\label{eq:diag-of-qSqt}
\sum_{n\geq1} \frac{s_{n-1,n}}{n!} t^n  =
\diag_{q,t}\biggl(\sum_{m\geq1, \ n\geq0} \frac{s_{m-1,n}}{n!} q^m t^n\biggr) =
\diag_{q,t}(q S(q, t)) ,
\end{equation}
and similarly
\begin{equation}\label{eq:diag-of-qHqt}
\sum_{n\geq1} s_{n-1,n} t^n  = \diag_{q,t}(q H(q, t)) ,
\end{equation}
where we have introduced the notation
\begin{gather*}
\diag_{q,t}\biggl(\sum_{i,j}c_{i,j}q^it^j\biggr) = \sum_i c_{i,i}t^i
\end{gather*}
for the diagonal of a bivariate series.

Our third proof, by combinatorial recurrences (Section~\ref{sec:graphs}), is an adaptation of a method
initially set up by Read~\cite{Read-1970-SUE} and developed explicitly for 3-regular graphs by Wormald~\cite{Wormald-1979-ELG}:
not only do we have to generalize it from 3-regular graphs to $(3,1)$-regular graphs,
but also in a way that edges can be counted together with vertices.
In a nutshell, the method consists in finding the right transformations of graphs,
first removing an edge or a vertex then fixing the local shape while diminishing the number of vertices,
so as to lead to fixpoint equations.
In doing so, auxiliary classes of graphs are needed:
in particular, a limited number of vertices of degree~2 has to be allowed (0, 1, or~2)
and a linear differential system satisfied by the generating functions of four interrelated graph classes
is obtained.
Differential elimination then leads to differential equations for $(3,1)$-regular graphs.

The three sections \ref{sec:triplesum}, \ref{sec:residues}, and~\ref{sec:graphs} are mutually independent.
Depending on their familiarity and mathematical preference,
the reader may decide to read
either the shorter algebraic approach (Section~\ref{sec:residues}) or the more combinatorial approach (Section~\ref{sec:graphs}).
Still, we will start with the more involved proof by a triple sum (Section~\ref{sec:triplesum})
as it appeared to be the method one would try by reflex when faced with the problem.
Finally, in Section~\ref{sec:conclusion} we compare the three proofs and their different relation to computation.

We supplement this text with various scripts as an online archive available at \url{https://files.inria.fr/chyzak/conj16/}.

\section{Proof by a triple-sum representation}
\label{sec:triplesum}

In this section we formulate a triple-sum representation for~$a_k$
before applying creative telescoping to produce the recurrence relation~\eqref{eq:conj-rec}.
Let us repeat that we could not derive any complete proof from the classical algorithms alone.

\subsection{Formulation as a sum}

We will obtain the numbers~$a_k$ from the even part of the series~\eqref{eq:diag-of-qHqt},
for which we first need to determine the ordinary generating function~$H(q, t)$.
In the spirit of Gessel's work~\cite[Theorem~7]{Gessel-1990-SFP},
we can view the representation~\eqref{eq:ord-gf-as-scalp} of~$H(q, t)$ by a scalar product
as the evaluation of a Hadamard product (with respect to the monomials in $p_1,p_2,p_3$),
\begin{equation*}
H(q,t) = \left(\exp(f(q, p)) \odot \frac{1}{1-tg(p)} \odot K(1p_1)K(2p_2)K(3p_3)\right)_{p_1=p_2=p_3=1}
\end{equation*}
where $K(x) = \sum_{\ell=0}^\infty \ell! x^\ell$,
or in a form that simplifies the calculation to come,
\begin{equation}\label{eq:H-as-Hadamard-products}
H(q,t) = \left(\exp(f(q, p)) \odot \frac{1}{1-tg(1p_1, 2p_2, 3p_3)} \odot K(p_1)K(p_2)K(p_3)\right)_{p_1=p_2=p_3=1} .
\end{equation}

In order to determine the coefficient $[q^m t^n] H(q, t)$ for related $m$ and~$n$,
we compute two extractions separately.
The following calculation is reminiscent of the calculation for 3-regular graphs
in \cite[(6.1), formula for~$C_n^*$]{Read-1960-ELR},
a variant of the original calculation for multigraphs in \cite[Section~7(e), obtaining~$K_n$]{Read-1959-ELR}.
The analogue of~$g$ for 3-regular graphs has three terms,
hence the double sums in those references,
while the polynomial~$g$ we use here has four terms,
leading to triple sums.

First, multiple applications of the binomial theorem yield:
\begin{multline*}
[t^n] \frac{1}{1-tg(1p_1, 2p_2, 3p_3)} = \left(p_1 + \tfrac16 p_1^3 + p_1p_2 + p_3\right)^n = \\
\sum_{a=0}^n\sum_{b=0}^{n-a}\sum_{c=0}^{n-a-b} \binom{n}{a,b,c,n-a-b-c}
  \frac{1}{2^b3^b} p_1^{a+3b+c} p_2^c p_3^{n-a-b-c} .
\end{multline*}

Second, we split $\exp(f(q, p))$ into a product of three exponentials
and extract the coefficient of~$q^{3u}$ in the last one,
next the coefficient of~$q^{2v}$ in the middle one,
so as to derive, using the binomial theorem again:
\begin{multline*}
[q^m] \exp(f(q, p)) = [q^m] \exp(\tfrac12 q (p_1^2 - p_2)) \exp(- \tfrac14 q^2 p_2^2) \exp(\tfrac16 q^3 p_3^2)= \\
\sum_{u=0}^{\lfloor m/3\rfloor} \sum_{v=0}^{\lfloor (m-3u)/2\rfloor}
  \frac{1}{(m-3u-2v)! \, v! \, u!}
  \left(\frac{p_1^2-p_2}{2}\right)^{m-3u-2v} \left(-\frac{p_2^2}{4}\right)^v \left(\frac{p_3^2}{6}\right)^u = \\
\sum_{u=0}^{\lfloor m/3\rfloor} \sum_{v=0}^{\lfloor (m-3u)/2\rfloor} \sum_{w=0}^{m-3u-2v}
  \frac{1}{u! \, v! \, w! \, (m-3u-2v-w)!}
  \frac{(-1)^{m+u+v+w}}{2^{m-2u} 3^u} \\
  {}\times p_1^{2w} p_2^{m-3u-w} p_3^{2u}
.
\end{multline*}
Note that, in the resulting triple sum,
$v$~appears only in the coefficients of the monomials in the~$p_i$,
not in the exponents of the~$p_i$.

At this point, in view of the diagonal~\eqref{eq:diag-of-qHqt} we are interested in
and of the representation~\eqref{eq:H-as-Hadamard-products},
we force $n=2k$ and~$m=2k-1$ for~$k\geq1$,
then, in order to determine terms contributing to the Hadamard product, we force
\begin{equation*}
p_1^{a+3b+c} p_2^c p_3^{2k-a-b-c} = p_1^{2w} p_2^{2k-1-3u-w} p_3^{2u} .
\end{equation*}
A calculation proves the equivalent conditions
\begin{equation}\label{eq:required-abc}
a = k + 1 , \qquad b = u + w - k , \qquad c = 2k - 1 -3u - w .
\end{equation}
So, owing to~\eqref{eq:H-as-Hadamard-products},
the coefficient of $q^{2k-1} t^{2k}$ in~\eqref{eq:diag-of-qHqt}, which is nothing but~$a_k$, is equal to
\begin{multline*}
\sum_{u=0}^{\lfloor (2k-1)/3\rfloor} \sum_{v=0}^{\lfloor (2k-1-3u)/2\rfloor} \sum_{w=0}^{2k-1-3u-2v}
  \binom{2k}{a,b,c,2k-a-b-c} \frac{1}{2^b 3^b} \\
\times
  \frac{1}{u! \, v! \, w! \, (2k-1-3u-2v-w)!} \frac{(-1)^{2k-1+u+v+w}}{2^{2k-1-2u} 3^u} \\
\times
  (2w)! \, (2k-1-3u-w)! \, (2u)!
\end{multline*}
where the final product of factorials reflects the second Hadamard product with functions~$K$ in~\eqref{eq:H-as-Hadamard-products}
and where $a,b,c$ are given by~\eqref{eq:required-abc}.
The number~$a_k$ is thus equal to
\begin{multline*}
\sum_{u=0}^{\lfloor (2k-1)/3\rfloor} \sum_{v=0}^{\lfloor (2k-1-3u)/2\rfloor} \sum_{w=0}^{2k-1-3u-2v}
  \frac{(-1)^{u+v+w+1}}{2^{k-1-u+w} 3^{2u+w-k}} \\
\times
  \binom{2k}{k+1,u+w-k,2k-1-3u-w,2u} \\
\times
  \frac{(2w)! \, (2k-1-3u-w)! \, (2u)!}{u! \, v! \, w! \, (2k-1-3u-2v-w)!}
.
\end{multline*}

We can simplify this sum by replacing the multinomial coefficient by factorials.
The summation ranges ensure that all of its lower arguments
are non-negative, except maybe~$u+w-k$.
But the constraint $0 \leq u \leq \lfloor (2k-1)/3\rfloor$ implies $1 \leq \lceil (k+1)/3\rceil \leq k-u  \leq k$,
so the sum with respect to~$w$ really starts at~$k-u$
and like the other, the argument~$u+w-k$ is non-negative.
Comparing the upper and the second lower argument of the multinomial, we get that
a non-zero contribution also requires~$u+w-k \leq 2k$,
or equivalently~$w \leq 3k-u$.
Observing $2k-1-3u-2v = (3k-u) - (k+1+2u+2v) < 3k-u$ we get that
the summation with respect to~$w$ can be reduced to the possibly-empty integer interval $\{k-u, \dots, 2k-1-3u-2v\}$,
and this yields
\begin{multline}\label{eq:last-not-natural}
a_k = \sum_{u=0}^{\lfloor (2k-1)/3\rfloor} \sum_{v=0}^{\lfloor (2k-1-3u)/2\rfloor} \sum_{w=k-u}^{2k-1-3u-2v}
  \frac{(-1)^{u+v+w+1}}{2^{k-1-u+w} 3^{2u+w-k}} \\
\times
  \frac{(2k)! \, (2w)!}{(k+1)! \, (u+w-k)! \, u! \, v! \, w! \, (2k-1-3u-2v-w)!}
\end{multline}
where the innermost sum is meant to be zero if~$k-u>2k-1-3u-2v$,
that is, when $2(u+v)\geq k$.

Note that all equivalent representations of the summand as products of binomials that we have tried
include remaining factorials or binomial coefficients in the denominator:
as such the summand does not lend itself to the method of binomial sums in the sense of \cite{BostanLairezSalvy-2013-CTR},
which otherwise would be the method of choice.

The next step of the creative-telescoping methodology is to derive a recurrence relation on the summand
that lends itself to an operation that transforms this relation into a recurrence relation for the triple sum.
To avoid the intractability of the many boundary terms and other exceptional terms that would be produced
from using~\eqref{eq:last-not-natural} directly
and from divisions by zero at singularities of the expressions,
both in~\eqref{eq:last-not-natural} and in intermediate calculations to come,
we replace (partially defined) factorials in the numerator and denominator
by using (total) functions defined by cases.
Namely,
we extend the factorial function to make it be zero on negative integers,
\begin{equation*}
n! = \begin{cases}
\Gamma(n+1) & \text{if $n \geq 0$}, \\
0 & \text{if $n < 0$}.
\end{cases}
\end{equation*}
and we introduce
\begin{equation}\label{eq:inverse-factorial}
\invfact{n} = \begin{cases}
1/n! & \text{if $n \geq 0$}, \\
0 & \text{if $n < 0$}.
\end{cases}
\end{equation}
We will call this last function the “inverse factorial”, although we can only derive the formula
\begin{equation*}
n! \times \invfact{n} = \begin{cases}
1 & \text{if $n \geq 0$}, \\
0 & \text{if $n < 0$},
\end{cases}
\end{equation*}
instead of an unconditional simplification to~1.
At least, the inverse factorial behaves better than the factorial:
it satisfies the recurrence
\begin{equation}\label{eq:rec-for-inverse-factorial}
\invfact{n} = (n+1) \invfact{(n+1)}
\end{equation}
for all~$n \in \bZ$,
while
\begin{equation}\label{eq:guarded-rec-for-factorial}
n! = n (n-1)!
\end{equation}
holds only for non-zero~$n \in \bZ$.

Our motivation to introduce the inverse factorial is
that \eqref{eq:last-not-natural}~can now be rewritten
as a finite sum “with natural bounds”,
that is, without the need to write explicit bounds, in the form
\begin{equation}\label{eq:triple-sum-with-natural-bounds}
a_k = \sum_{u\in\bZ} \sum_{v\in\bZ} \sum_{w\in\bZ} f_{k,u,v,w}
\end{equation}
where
\begin{multline}\label{eq:def-f}
f_{k,u,v,w} =
  (-1)^{u+v+w+1} 2^{1+u-w-k} 3^{k-2u-w} \\
\times
  (2k)! \, (2w)! \times \invfact{(k+1)} \, \invfact{(u+w-k)} \, \invfact{u} \, \invfact{v} \, \invfact{w} \, \invfact{(2k-1-3u-2v-w)}
.
\end{multline}

\subsection{Shortcomings of conventional creative telescoping}
\label{sec:shortcomings}

Existing creative-telescoping algorithms for sums only deal with single sums,
with two exceptions:
for double sums of hypergeometric terms~\cite{ChenHouMu-2006-TMD}
and for triple sums of rational terms~\cite{ChenDuFang-2025-SSM}.
So for triple hypergeometric sums, one has either to appeal
to an iterated use of Chyzak's algorithm~\cite{Chyzak-2000-EZF},
or to use Koutschan's ansatz~\cite{Koutschan-2010-FAC},
which is often much faster but not proved to be complete,
even on the class of holonomic hypergeometric terms.

Whether we run Chyzak's Maple implementation \texttt{Mgfun}%
\footnote{The package is available at \url{https://mathexp.eu/chyzak/mgfun.html}. Use the command \texttt{creative\_telescoping}.}
or Koutschan's Mathematica implementation \texttt{HolonomicFunctions}%
\footnote{The package is available at \url{https://risc.jku.at/sw/holonomicfunctions/}. Use the command \texttt{CreativeTelescoping} three times.}
of Chyzak's creative-telescoping algorithm of~\cite{Chyzak-2000-EZF},
we obtain a set of seven pairs of operators that represent linear recurrence relations between values of the sequence~$(f_{k,u,v,w})_{k,u,v,w}$
and values of the inner single and double sums in~\eqref{eq:triple-sum-with-natural-bounds}.
Those relations are parametrized by~$k,u,v,w$ and ``promised'' to be valid for ``most'' integer values of the parameters.
There are two obstacles to use the obtained recurrence relations in a formally complete proof:
first, the algorithm does not output the domain of validity of the equations;
and second, using the relations on the sequence to get relations on the triple sum is,
to the best of our effort, intractable.

To further comment on this,
write $\cR_{k,u,v,w}$ for the algebra generated over $\bQ(k,u,v,w)$ by symbols $S_k,S_u,S_v,S_w$ satisfying
\[ S_k k = (k+1) S_k, \quad S_u u = (u+1) S_u, \quad S_v v = (v+1) S_v, \quad S_w w = (w+1) S_w,\]
all other pairs of symbols taken from $k,u,v,w,S_k,S_u,S_v,S_w$ commuting.
This algebra is an Ore algebra in the sense of the original definition in~\cite{ChyzakSalvy-1998-NCE},
which highlights the commutativity of the symbols $S_k,S_u,S_v,S_w$
to make a non-commutative analogue of the Gröbner-basis theory available.
We also consider $\cR_{k,u,v}$ for the algebra generated over~$\bQ(k,u,v)$ by symbols $S_k,S_u,S_v$ satisfying the same relations,
and similarly defined $\cR_{k,u}$ and~$\cR_k$ as well.
As usual, elements of~$\cR_{k,u,v,w}$ are operators that act on an abstract difference module,
but not on genuine sequences.
In our case, the sequence~$f$ corresponds to an element~$\hat f$
in a $\bQ(k,u,v,w)$-vector space endowed with an action of~$\cR_{k,u,v,w}$
so that $\hat f$~satisfies \emph{absolutely} the first-order recurrence relations reflected by the operators
\begin{align}
\label{eq:S_k-C_1}
  Z_1 := S_k &- C_1  &  \text{for } C_1 &= \frac{3(k+1)(2k+1)(u+w-k)}{(k+2)\prod_{i=0}^1(2k+i-3u-2v-w)} , \\
\label{eq:S_u-C_2}
  Z_2 := S_u &- C_2  &  \text{for } C_2 &= - \frac{2\prod_{i=1}^3(2k-i-3u-2v-w)}{9(u+w+1-k)(u+1)} , \\
\label{eq:S_v-C_3}
  Z_3 := S_v &- C_3  &  \text{for } C_3 &= - \frac{\prod_{i=1}^2(2k-i-3u-2v-w)}{v+1} , \\
\label{eq:S_w-C_4}
  Z_4 := S_w &- C_4  &  \text{for } C_4 &= - \frac{(2w+1)(2k-1-3u-2v-w)}{3(u+w+1-k)} .
\end{align}%
\noeqref{eq:S_u-C_2,eq:S_v-C_3}%
In contrast, the sequence~$f$ only satisfies the same recurrence relations where denominators do not vanish,
for example, $f_{k,u+1,v,w}$~is $C_2(k,u,v,w) f_{k,u,v,w}$ only if $k\neq u+w+1$.
Then, calculating with either implementation, the set of seven obtained pairs decomposes as follows:
\begin{itemize}
\item for three pairs $(P_i,Q_i)$, $i=1,2,3$,
  the combination $P_i - \Delta_w Q_i$ is an element of the left ideal~$\hat I$ in~$\cR_{k,u,v,w}$
  generated by the recurrence operators \eqref{eq:S_k-C_1}--\eqref{eq:S_w-C_4};
\item for three pairs $(P'_i,Q'_i)$, $i=1,2,3$,
  the combination $P'_i - \Delta_v Q'_i$ is an element of the left ideal~$I$ in~$\cR_{k,u,v}$ generated by $P_1,P_2,P_3$;
\item for a single pair $(P'',Q'')$,
  the combination $P'' - \Delta_u Q''$ is an element of the left ideal~$I'$ in~$\cR_{k,u}$ generated by $P'_1,P'_2,P'_3$.
\end{itemize}
These pairs contain Gröbner bases for the ideals mentioned above:
$\{P_1,P_2,P_3\}$ for~$\hat I$,
$\{P'_1,P'_2,P'_3\}$ for~$I$,
and $\{P''\}$~for~$I'$.
From this data, we obtain an explicit and direct representation of~$P''$ as follows.
First, we have:
\begin{align*}
P_i - \Delta_w Q_i &= \sum_{\ell=1}^4 U_{i,\ell} Z_\ell, \qquad (i=1,2,3) , \\
P'_j - \Delta_v Q'_j &= \sum_{i=1}^3 U'_{j,i} P_i, \qquad (j=1,2,3) , \\
P'' - \Delta_u Q'' &= \sum_{j=1}^3 U''_j P'_j , \\
\end{align*}
where $Z_1,Z_2,Z_3,Z_4$ are the first-order elements in~$\cR_{k,u,v,w}$ mentioned in the first item
and where
\begin{equation*}
U_{i,\ell} \in \cR_{k,u,v,w} , \quad
U'_{j,i} \in \cR_{k,u,v} , \quad
U''_j \in \cR_{k,u} .
\end{equation*}
Second, eliminating the $P_i$ and~$P'_j$ from the equations above results in
\begin{multline}\label{eq:rec-as-deltas}
P'' = \Delta_u Q''
    + \Delta_v \Biggl(\sum_{1\leq j\leq3} U''_j Q'_j\Biggr)
    + \Delta_w \Biggl(\sum_{1\leq i,j\leq3} U''_j U'_{j,i} Q_i\Biggr) \\
    \qquad\qquad\qquad\qquad\qquad\qquad\qquad
     {} + \sum_{\ell=1}^4 \Biggl(\sum_{1\leq i,j\leq3} U''_j U'_{j,i} U_{i,\ell}\Biggr) Z_\ell \\
\shoveleft{\phantom{P''}
    = \Delta_u Q'' + \Delta_v Q' + \Delta_w Q + \sum_{\ell=1}^4 C_\ell Z_\ell ,}
    \qquad\qquad\qquad\qquad\qquad\qquad\qquad
\end{multline}
in which each of $Q', Q, C_\ell$ is defined to be equal to the parenthesized term with the same role in the previous line.

The method of creative telescoping next amounts
to applying~\eqref{eq:rec-as-deltas} to the sequence~$(f_{k,u,v,w})_{k,u,v,w}$
before summing over all relative integers $u,v,w$,
in the hope that the resulting right-hand side is zero for all~$k \in \bZ_{\geq0}$.
The algorithm has proved that
applying~\eqref{eq:rec-as-deltas} to the term~$\hat f$ is zero,
and this term~$\hat f$ satisfies more recurrence relations than the sequence,
or recurrence relations without validity constraint,
so some proof is needed.
A problem is that the recurrence relation
\begin{equation}\label{eq:conditional-rec}
\begin{aligned}
(P''\cdot f)(k,u,v,w)
  &= (Q''\cdot f)(k,u+1,v,w) - (Q''\cdot f)(k,u,v,w) \\
  &+ (Q'\cdot f)(k,u,v+1,w) - (Q'\cdot f)(k,u,v,w) \\
  &+ (Q\cdot f)(k,u,v,w+1) - (Q\cdot f)(k,u,v,w) \\
  &+ \sum_{\ell=1}^4 (C_\ell Z_\ell\cdot f)(k,u,v,w)
\end{aligned}
\end{equation}
does not in general hold for all $(k,u,v,w) \in \bZ^4$:
while the left-hand side is always well-defined by construction,
$Q''$, $Q'$, $Q$ and the~$C_\ell$ may well involve denominators that make the evaluations meaningless,
and additionally, even if no~$Z_\ell$ had a non-trivial denominator,
the identity $(Z_\ell\cdot f)(k,u,v,w)$ would not need to hold globally,
owing to the proviso~$k\neq0$ on the recurrence relation~\eqref{eq:guarded-rec-for-factorial} satisfied by the factorial sequence.
Because of this
we need to remove some points of the summation ranges from the set at which we specialize~\eqref{eq:conditional-rec}
before summation,
and these points occur in intersecting infinite families.
To proceed,
we would need to decompose the expression into a number of single and double sums,
according to the geometry of the resulting ranges.
In contrast with our original sum~\eqref{eq:triple-sum-with-natural-bounds},
which we managed to formulate with natural bounds,
the many single and double sums are over ranges with boundaries.
The whole geometry is so complicated that we have not been able to make it explicit,
and we declared that the resulting expression is unmanageable.

The temporary conclusion of this section is that, despite our many efforts,
we have not been able to prove the conjecture by applying creative telescoping to the triple sum~\eqref{eq:last-not-natural}.

\subsection{Completing a proof}
\label{sec:triplesum-ok}

To avoid all sorts of difficulties or gaps in our proof by creative telescoping for multiple sums,
we adapt the method to compute a different form of a creative-telescoping relation,
improving on~\eqref{eq:conditional-rec},
from which a formally complete proof will be possible.
We proceed in several steps.

\paragraph{Initial setup.}

We introduce the left ideal~$\hat I$ generated in~$\cR_{k,u,v,w}$ by \eqref{eq:S_k-C_1}--\eqref{eq:S_w-C_4}.
Computing modulo those operators disregarding any domain of validity
is what computer-algebra systems do when simplifying expression representing sequences%
\footnote{For example, the Maple command \texttt{simplify(binomial(n+1,k) / binomial(n,k))} outputs $(n+1)/(n+1-k)$.}%
.
Indeed, the operators \eqref{eq:S_k-C_1}--\eqref{eq:S_w-C_4}
do not represent recurrence relations on~$f_{k,u,v,w}$
that would be valid for all integer $k$, $u$, $v$, and~$w$,
as, for example, the relation
\begin{equation*}
  f_{k+1,u,v,w} = C_1(k,u,v,w) f_{k,u,v,w}
\end{equation*}
holds only when the denominator of~$C_1$ is non-zero.
We could prove that the variant relation obtained after clearing denominators
is valid for all integer~$(k,u,v,w)$, and similarly for the three other operators,
but the fact is that the next calculation will use
rational multiples of the above and introduce further denominators anyway,
owing to purely algebraic simplifications,
so we do not follow this route.
Also, the proof technique that we would use is an instance of the method
that we will develop in more generality below.

\paragraph{Naive difference-algebra ansatz.}

By using the \texttt{FindCreativeTelescoping} command in Koutschan's \texttt{HolonomicFunctions} package,
we compute 4-variate rational functions $Q_1$, $Q_2$, and~$Q_3$ such that
\begin{equation*}
  P(k,S_k) - \Delta_u Q_1(k,u,v,w) - \Delta_v Q_2(k,u,v,w) - \Delta_w Q_3(k,u,v,w) \in \hat I
\end{equation*}
for the operator $P = p_5(k)S_k^5+\dots+p_0(k)$ underlying the left-hand side of~\eqref{eq:conj-rec}.
At this point we would like to apply this operator to~$f_{k,u,v,w}$ and sum over~$(u,v,w)$ in~$\bZ^3$,
but the rational functions have several factors in their denominators that make the evaluation meaningless at many integer choices for~$(k,u,v,w)$.
For example, $Q_1$~has the denominator
\begin{equation*}
  (3k+76)(k+3)\prod_{i=2}^8 (2k+i-3u-2v-w) ,
\end{equation*}
which beside~$k=-3$ has for large enough~$k$ a number of integer zeros proportional to~$k^2$.

\paragraph{Removal of denominators.}

We modify the previous rational functions to obtain new rational functions $R_1$, $R_2$, and~$R_3$
satisfying
\begin{multline*}
  P(k,S_k) - \Delta_u R_1(k,u,v,w) S_k^9S_w^9 - \Delta_v R_2(k,u,v,w) S_k^{10}S_w^{10} \\
  - \Delta_w R_3(k,u,v,w) S_k^{10}S_w^{10} \in \hat I ,
\end{multline*}
with the purpose that the denominators of the~$R_i$ have as few integer zeros as possible.
This is done by forcing $R_1(k,u,v,w) S_k^9S_w^9 - Q_1(k,u,v,w) \in \hat I$, or equivalently
\begin{equation*}
  R_1(k,u,v,w) \biggl(\prod_{i=0}^8 C_4(k+9,u,v,w+i)\biggr) \biggl(\prod_{i=0}^8 C_1(k+i,u,v,w)\biggr) = Q_1(k,u,v,w) ,
\end{equation*}
and similar relations for $R_2$ and~$R_3$.
The only integer zeros of any of the~$R_i$ are for tuples $(k,u,v,w)$ with~$k=-3$ and unrestricted integers $u$, $v$, and~$w$,
a situation that we forged by trial and error
by testing monomials in the shifts that ended up being chosen to be $S_k^9S_w^9$ and~$S_k^{10}S_w^{10}$.
So at this point, the expression
\begin{multline}\label{eq:to-be-zero}
  z_{k,u,v,w} :=
  p_5(k) f_{k+5,u,v,w} + p_4(k) f_{k+4,u,v,w} + p_3(k) f_{k+3,u,v,w} \\
  + p_2(k) f_{k+2,u,v,w} + p_1(k) f_{k+1,u,v,w} + p_0(k) f_{k,u,v,w} \\
  - R_1(k,u+1,v,w) f_{k+9,u+1,v,w+9} + R_1(k,u,v,w) f_{k+9,u,v,w+9} \\
  - R_2(k,u,v+1,w) f_{k+10,u,v+1,w+10} + R_2(k,u,v,w) f_{k+10,u,v,w+10} \\
  - R_3(k,u,v,w+1) f_{k+10,u,v,w+11} + R_3(k,u,v,w) f_{k+10,u,v,w+10}
\end{multline}
is well defined on any integer $k$, $u$, $v$, and~$w$ unless~$k=-3$.
We therefore consider proving it is zero under the sole restriction~$k\neq-3$.

\paragraph{Rewriting in terms of a common shift.}

Writing each~$f_{k+p,u+q,v+r,w+s}$ appearing in~\eqref{eq:to-be-zero}
in the form of a rational function times~$f_{k,u,v,w}$
would re-introduce denominators and make the evaluation of~$z_{k,u,v,w}$
undefined at many integer choices of~$(k,u,v,w)$.
The way to avoid this, which we determined by trial and error, is
not only to rewrite them in terms of~$f_{k+11,u,v,w+12}$ instead of~$f_{k,u,v,w}$,
but also to consider relations without denominators.
For example, from
\begin{equation*}
S_k^{10}S_w^{11} + \frac{(k+12)(2k+9-3u-2v-w)}{(k+11)(2k+21)(2w+23)} S_k^{11}S_w^{12} \in \hat I
\end{equation*}
which reflects a simplification of binomial expressions in a computer-algebra system
and can be otherwise derived from the elements \eqref{eq:S_k-C_1}--\eqref{eq:S_w-C_4} of~$\hat I$,
we expect the relation
\begin{multline}\label{eq:one-case-of-rewrite}
(k+11)(2k+21)(2w+23) f_{k+10,u,v,w+11} = \\
  -(k+12)(2k+9-3u-2v-w) f_{k+11,u,v,w+12} .
\end{multline}
Analyzing the domain of validity of this relation and of the eleven other relations
that will play a similar role to simplify~\eqref{eq:to-be-zero}
is the bottleneck of our process.
We developed a computational approach whose description we postpone to the end of the current section:
the bottomline is that
each of the twelve terms of the form~$f_{k+p,u+q,v+r,w+s}$ that appears in~\eqref{eq:to-be-zero}
is related to~$f_{k+11,u,v,w+12}$ by a relation similar to~\eqref{eq:one-case-of-rewrite}
that holds for all integer $k$, $u$, $v$, and~$w$.

\paragraph{Reducing the candidate relation to zero.}

We process the terms~$f_{k+p,u+q,v+r,w+s}$ independently.
As an example,
the numerator of the coefficient~$R_3(k,u,v,w+1)$ of~$f_{k+10,u,v,w+11}$ in~$z_{k,u,v,w}$ is
\begin{equation*}
(k+11)(2k-1-3u-2v-w)(2k-3u-2v-w)\times\text{(large irreducible polynomial)} .
\end{equation*}
Note the presence of~$k+11$, but not of $(2k+21)(2w+23)$, so that
after multiplying~\eqref{eq:to-be-zero} by the latter,
the term involving~$f_{k+10,u,v,w+11}$ can be rewritten
using~\eqref{eq:one-case-of-rewrite} without introducing any additional denominators.
As a whole,
proceeding in the same manner with all terms~$f_{k+p,u+q,v+r,w+s}$ appearing in~\eqref{eq:to-be-zero}
requires pre-multiplying~$z_{k,u,v,w}$ by the factor
\begin{equation*}
p(k,w) := (k+3) \biggl(\prod_{i=5}^{10} (2k+2i+1)\biggr) (3k+76) \biggl(\prod_{i=0}^{11} (2w+2i+1\biggr) .
\end{equation*}
Doing so, we rewrite~$p(k,w) z_{k,u,v,w}$ into the form~$E(k,u,v,w) f_{k+11,u,v,w+12}$
where the coefficient~$E$ is some rational expression in the variables~$k,u,v,w$.
Observe that $p(k,w)$~is non-zero for any integer choice of~$(k,u,v,w)$ provided~$k\neq-3$.
So, to prove~$z_{k,u,v,w} = 0$ for~$k\neq-3$, it is sufficient
to prove~$p(k,w) z_{k,u,v,w} = 0$ for all integer $k$, $u$, $v$, and~$w$.
Normalizing the rational function~$E$, we obtain zero, thus proving this latter goal,
and as, of course, the normalization introduces no additional denominators,
the telescoping relation~$z_{k,u,v,w} = 0$
holds for all integer $k$, $u$, $v$, and~$w$ provided~$k\neq-3$.

\paragraph{Telescoping the telescoping relation.}

Fixing~$k\neq-3$ and summing~\eqref{eq:to-be-zero}, which has been proved to be zero,
over all integer $u$, $v$, and~$w$,
yields~\eqref{eq:conj-rec},
where $a_k$~is the (actually finite) sum defined by~\eqref{eq:triple-sum-with-natural-bounds}.

\paragraph{Proving two-term relations between the $f_{k+p,u+q,v+r,w+s}$ and~$f_{k+11,u,v,w+12}$.}

We automated the proof of the relation~\eqref{eq:one-case-of-rewrite} and its siblings
for all integer $k$, $u$, $v$, and~$w$.
To this end, we first considered developing a data structure for guarded expressions,
to follow ideas from~\cite{DolzmannSturm-1997-GEP},
but this proved to be not so convenient.
This convinced us that we should use unguarded rewriting rules,
valid for all integer $k$, $u$, $v$, and~$w$.
To this end, we first generalized the guarded~\eqref{eq:guarded-rec-for-factorial} in a form
that arithmetizes its piecewise nature:
using Iverson's bracket notation,
\begin{equation*}
\iverson{\cP} = \begin{cases}
1 & \text{if the predicate $\cP$ holds}, \\
0 & \text{if it does not hold},
\end{cases}
\end{equation*}
we obtain the relation
\begin{equation}\label{eq:unguarded-rec-for-factorial}
n! = n (n-1)! + \iversonzero{n}
\end{equation}
that holds for all~$n \in \bZ$.
By the nature of our calculations, Iverson brackets will only have arguments that are affine equations
in $k$, $u$, $v$, and~$w$ with integer coefficients.
A nice property of this is that expressions can be simplified by the formula
\begin{equation}\label{eq:simplify-Iverson}
c(e) \, \iversonzero{e} = c(0) \, \iversonzero{e} ,
\end{equation}
valid when $e$~takes values from~$\bZ$.
In practice, we will make a choice of a symbol occurring in~$e$
and make a substitution that eliminates it from~$c(e)$.
For example, $c(k-w) \, \iversonzero{k-u-v}$ simplifies to $c(u+v-w) \, \iversonzero{k-u-v}$.
To explain the nature of our automated proof, we exemplify it here
by proving~\eqref{eq:one-case-of-rewrite} for all integer $k$, $u$, $v$, and~$w$
by hand, through the following calculation:
\begin{align*}
-(k+12)&(2k+9-3u-2v-w) f_{k+11,u,v,w+12} \mystackrel1= \\
-(k&+12)(2k+9-3u-2v-w) \\
  &\times (-1)^{u+v+w+1} 2^{-22+u-w-k} 3^{k-1-2u-w} \\
  &\times (2k+22)! \, (2w+24)! \\
  &\times \invfact{(k+12)} \, \invfact{(u+w+1-k)} \, \invfact{u} \, \invfact{v} \, \invfact{(w+12)} \, \invfact{(2k+9-3u-2v-w)} \mystackrel2= \\
(-1)&^{u+v+w} 2^{-22+u-w-k} 3^{k-1-2u-w} \\
  &\times (2k+22)! \, (2w+24)! \\
  &\times \invfact{(k+11)} \, \invfact{(u+w+1-k)} \, \invfact{u} \, \invfact{v} \, \invfact{(w+12)} \, \invfact{(2k+8-3u-2v-w)} \mystackrel3= \\
(-1)&^{u+v+w} 2^{-22+u-w-k} 3^{k-1-2u-w} \\
  &\times \bigl(2(k+11)(2k+21) \, (2k+20)! + \iversonzero{k+11}\bigr) \\
  &\times \bigl(2(w+12)(2w+23) \, (2w+22)! + \iversonzero{w+12}\bigr) \\
  &\times \invfact{(k+11)} \, \invfact{(u+w+1-k)} \, \invfact{u} \, \invfact{v} \, \invfact{(w+12)} \, \invfact{(2k+8-3u-2v-w)} \mystackrel4= \\
(-1)&^{u+v+w} 2^{-20+u-w-k} 3^{k-1-2u-w} \\
  &\times (k+11)(2k+21)(w+12)(2w+23) \, (2k+20)! \, (2w+22)! \\
  &\times \invfact{(k+11)} \, \invfact{(u+w+1-k)} \, \invfact{u} \, \invfact{v} \, \invfact{(w+12)} \, \invfact{(2k+8-3u-2v-w)} \\
+ (-&1)^{u+v+w} 2^{-11+u-w} 3^{-12-2u-w} \\
  &\times \iversonzero{k+11} \, 2(w+12)(2w+23) \, (2w+22)! \\
  &\times \invfact{0} \, \invfact{(u+w+12)} \, \invfact{u} \, \invfact{v} \, \invfact{(w+12)} \, \invfact{(-14-3u-2v-w)} \\
+ (-&1)^{u+v} 2^{-10+u-k} 3^{k+11-2u} \\
  &\times 2(k+11)(2k+21) \, (2k+20)! \, \iversonzero{w+12} \\
  &\times \invfact{(k+11)} \, \invfact{(u-11-k)} \, \invfact{u} \, \invfact{v} \, \invfact{0} \, \invfact{(2k+20-3u-2v)} \\
+ (-&1)^{u+v} 2^{1+u} 3^{-2u} \, \iversonzero{k+11} \, \iversonzero{w+12} \, \invfact{0} \, \invfact{u} \, \invfact{u} \, \invfact{v} \, \invfact{0} \, \invfact{(-2-3u-2v)} \mystackrel5= \\
(k+11)&(2k+21)(2w+23) \\
&\times (-1)^{u+v+w} 2^{-20+u-w-k} 3^{k-1-2u-w} \\
&\times (2k+20)! \, (2w+22)! \\
&\times \invfact{(k+11)} \, \invfact{(u+w+1-k)} \, \invfact{u} \, \invfact{v} \, \invfact{(w+11)} \, \invfact{(2k+8-3u-2v-w)} \mystackrel6= \\
(k+11)&(2k+21)(2w+23) f_{k+10,u,v,w+11} ,
\end{align*}
where
$\mystackrel1=$ and~$\mystackrel6=$ are by expanding the definition~\eqref{eq:def-f} of~$f$,
$\mystackrel2=$~is by applying~\eqref{eq:rec-for-inverse-factorial},
$\mystackrel3=$~is by applying~\eqref{eq:unguarded-rec-for-factorial},
$\mystackrel4=$~is by expanding into four terms the product of two factors involving Iverson brackets
before applying~\eqref{eq:simplify-Iverson},
and $\mystackrel5=$~is by realizing that three out of the four added terms are zero:
for example, $\invfact{(k+11)} \, \invfact{(u-11-k)} \, \invfact{u} \, \invfact{v} \, \invfact{0} \, \invfact{(2k+20-3u-2v)}$ being non-zero requires
$k+11\geq0$, $u-11-k\geq0$, $u\geq0$, $v\geq0$, $2k+20-3u-2v\geq0$,
and therefore
$0 \leq 2k+20-3u-2v = -2-2(u-11-k) -u-2v \leq -2$,
a contradiction;
the other two are proved similarly.

To deal with the many cases to be considered in the full proof that~$z_{k,u,v,w} = 0$,
we automated the simplification above in a Maple script
by developing a normalization procedure that:
\begin{itemize}
\item manipulates inputs that are polynomial expressions involving terms of the form $e!$, $\invfact{e}$, and~$\iversonzero{e}$ for various~$e$, and keeps intermediate polynomial expressions in a distributed collected representation with respect to those terms;
\item allows normalizations of rational functions in~$(k,u,v,w)$ only after verifying by inspection that no denominator has any integer zero;
\item for any term~$e!$ that appears in the input, tries to reduce the largest integer~$m$ such that the term~$(e+m)!$ also appears, by applying~\eqref{eq:unguarded-rec-for-factorial} with~$n=e+m-1$,
  and proceeds similarly with inverse factorials, by applying~\eqref{eq:rec-for-inverse-factorial};
\item inspects products of factorials and inverse factorials to remove products that are zero owing to a contradiction;
\item replaces a~$\Delta$ with~$0$ if its argument cannot be zero, as is for example the case with $\iversonzero{2k-3u+1}$;
\item systematically simplifies coefficients in front of a bracket~$\iversonzero{e}$ by~\eqref{eq:simplify-Iverson}.
\end{itemize}
To avoid any controlled use of a wrong rule,
we have used custom implementations of the functions $e!$, $\invfact{e}$, and~$\iversonzero{e}$.
To test for a contradiction in a product of factorials and inverse factorials,
we have encoded the existence of consistent arguments
as integer linear-programming problems.
For example, proving $\invfact{(k+11)} \, \invfact{(u-11-k)} \, \invfact{u} \, \invfact{v} \, \invfact{0} \, \invfact{(2k+20-3u-2v)} = 0$ boils down
to proving that the constraints $k+11\geq0$, $u-11-k\geq0$, $u\geq0$, $v\geq0$, $2k+20-3u-2v\geq0$
cannot be satisfied by integer specializations of the variables,
which algorithmically is solved by optimization of an arbitrary function in the corresponding polyhedron%
\footnote{In Maple, this can be done
by calling \texttt{Optimization:-LPSolve(0, \{k+11>=0, u-11-k>=0, u>=0, v>=0, 2*k+20-3*u-2*v>=0\}, assume=integer)},
which finds no feasible point.
Alternatively, one can call Maple's command \texttt{SMTLIB[Satisfiable]} with the option \texttt{logic="QF\_LIA"},
which calls the external solver~Z3.
However, this proved to be less efficient on our set of problems.}%
.
The resulting script for a complete proof is available in directory \texttt{triplesum/} of the online appendix.

\section{Proof by a residue representation}
\label{sec:residues}

To get another proof of~\eqref{eq:conj-rec},
we will look for a linear ODE with polynomial coefficients
for the exponential generating function
\begin{equation*}
A(t) = \sum_{k\geq1} \frac{a_k}{(2k)!} t^k ,
\end{equation*}
which is the even part of the diagonal~\eqref{eq:diag-of-qSqt} of the graph-counting bivariate series~$S(q, t)$ given by~\eqref{eq:Sqt-and-Hqt}.
That is, the exponential generating function~$A(t)$ satisfies
\begin{equation*}
A(t) = \even_t(\diag_{q,t}(q S(q, t))) ,
\end{equation*}
where we have introduced the notation
\begin{equation*}
\even_t\biggl(\sum_ic_it^i\biggr) = \sum_i c_{2i}t^i
\end{equation*}
for the even part of a univariate series.

An algorithm was provided in \cite{ChyzakMishna-2025-DES}
to compute a linear ODE satisfied by the univariate exponential generating function~$S(t) = S(1, t)$
counting vertex-labelled $(3,1)$-regular graphs with vertices marked by~$t$.
This algorithm is based on the already mentioned formula~\eqref{eq:St}.
It proceeds by reduction-based creative telescoping,
formulating an ad~hoc reduction procedure for symmetric scalar products.

A possible computation of an ODE for~$A(t)$ could therefore start
by modifying the procedure in~\cite{ChyzakMishna-2025-DES} so as to output a bivariate D-finite representation of~$S(q, t)$,
which would then have to be fed to an algorithm to compute an ODE for the diagonal,
before finally one would extract an ODE for the even part.
For future reference, let us mention that the differential system
is formed by the two bivariate linear differential operators
\begin{align}
\label{eq:bivariate-sys-only-dt}
&9q^3t^3(q^6t^4+2q^3t^2-2q^2t^2-2) \partial_t^2 \\
\notag
&\quad {} + (3q^{15}t^{10}+18q^{12}t^8-24q^{11}t^8+9q^9t^6-60q^8t^6+36q^7t^6 \\
\notag
&\qquad\qquad {} - 18q^6t^4+18q^5t^4-78q^3t^2+24q^2t^2+24) \partial_t \\
\notag
&\quad {} - (q^{17}t^{10}+4q^{14}t^8-8q^{13}t^8-28q^{10}t^6+36q^9t^6-8q^8t^6-8q^8t^4+40q^7t^4 \\
\notag
&\qquad\qquad {} + 36q^6t^4-64q^5t^4+16q^4t^4+4q^5t^2-96q^3t^2+40q^2t^2+24) \partial_t , \\
\label{eq:bivariate-sys-mixed}
&{} -2(q^6t^4+2q^3t^2-2q^2t^2-2)q \partial_q \\
\notag
&\qquad\qquad {} + 3(q^6t^4+2q^3t^2-2) t\partial_t - 4qt^2(q^3t^2-1) ,
\end{align}
where $\partial_q = d/dq$ and~$\partial_t = d/dt$,
and that the suggested diagonal computation returns
the very same linear differential operator~\eqref{eq:lodo-after-diag}
as the one obtained by the method we detail in this section.
A script for this calculation is available in directory \texttt{scalpbyreductions/} of the online appendix.

In a context where the edge count is not tracked,
another calculation of the linear ODE for~$S(t)$ was proposed in~\cite{BrochetChyzakLairez-2025-FMI},
especially in Section~6 of that reference.
After reinterpreting the scalar product~\eqref{eq:St} as a formal residue by the 3-fold residue formula
\begin{equation}\label{eq:St-as-res}
S(t) = \res_p(\exp(f(p)) \cL(\exp(t \tilde g(p)))) ,
\end{equation}
where
\begin{equation*}
\tilde g(p_1, p_2, p_3) = g(1 p_1, 2 p_2, 3 p_3)
\end{equation*}
and where $\cL({\cdot})$~denotes formal Laplace transform,
\begin{equation*}
\cL\biggl(\sum_{i,j,k,m,n}c_{i,j,k,m,n}p_1^ip_2^jp_3^kq^mt^n\biggr) = \sum_{i,j,k,m,n}c_{i,j,k,m,n}\frac{q^mt^n}{p_1^{i+1}p_2^{j+1}p_3^{k+1}} ,
\end{equation*}
a general-purpose creative-telescoping algorithm for computing integrals of holonomic modules is used.
This accommodates calculations of multiple integrals and multiple residues.
Instead of the previously proposed approach by generalizing~\cite{ChyzakMishna-2025-DES},
in which separating the diagonal step from the scalar product step looks like computing integrals iteratively,
the diagonal can be encoded as another residue.
Generalizing~\eqref{eq:St-as-res} in the form
\begin{equation*}
S(q, t) = \res_p(\exp(f(q, p)) \cL(\exp(t \tilde g(p))))
\end{equation*}
yields the 4-fold residue formula
\begin{equation*}
A(t) = \even_t(\res_q(S(q, q^{-1} t)))
= \even_t(\res_{p,q}(F \cL(G))) ,
\end{equation*}
where we have set
\begin{equation*}
F = \exp(f(q, p)) , \qquad G = \exp(q^{-1} t \tilde g(p)) .
\end{equation*}

From this point on, everything is computational,
after we derive a holonomic system of differential equations annihilating~$F \cL(G)$.
We proceed along the lines of~\cite{BrochetChyzakLairez-2025-FMI},
only adding simple modifications to take~$q$ into account.
As in that reference, the goal is to obtain differential skew polynomials
that generate an ideal in $\cD = \bQ(t)[q,p_1,p_2,p_3]\langle\partial_t,\partial_q,\partial_1,\partial_2,\partial_3\rangle$.

First, observe that $G$~is cancelled by the five skew polynomials
\begin{gather*}
q \partial_i - t \tilde g_i(p_1,p_2,p_3) , \ i=1,2,3 ,
\qquad
q^2 \partial_q + t \tilde g(p_1,p_2,p_3) ,
\qquad
q \partial_t - \tilde g(p_1,p_2,p_3)
\end{gather*}
where $\tilde g_i$~is the derivative of~$\tilde g$ with respect to~$p_i$.
Next, as detailed in~\cite{BrochetChyzakLairez-2025-FMI},
the formal Laplace transform~$\cL$ maps $p_i$~to~$\partial_i$ and $\partial_i$~to~$-p_i$,
in the sense that $\cL(G)$~is cancelled up to series with zero residue by the five skew polynomials
\begin{gather*}
q p_i - t \tilde g_i(-\partial_1,-\partial_2,-\partial_3) , \ i=1,2,3 ,
\\
q^2 \partial_q + t \tilde g(-\partial_1,-\partial_2,-\partial_3) ,
\qquad
q \partial_t - \tilde g(-\partial_1,-\partial_2,-\partial_3) .
\end{gather*}
Here and later on, $\tilde g(L_1,L_2,L_3)$~denotes the result of substituting
$L_1$ for~$p_1$, $L_2$ for~$p_2$, $L_3$ for~$p_3$,
all three being skew polynomials from~$\cD$,
in the commutative polynomial~$\tilde g(p_1,p_2,p_3)$.
Further, again as detailed in~\cite{BrochetChyzakLairez-2025-FMI},
any equation $P(\partial_1,\partial_2,\partial_3,\partial_q)\cdot H = 0$
leads after multiplication by~$F$ to
$P(\partial_1-f_1,\partial_2-f_2,\partial_3-f_3,\partial_q-f_q)\cdot (FH) = 0$,
where $f_i$~is the derivative of~$f$ with respect to~$p_i$ and $f_q$~that with respect to~$q$.
Thus, $F\cL(G)$~is cancelled up to series with zero residue by all of
\begin{gather}
\label{eq:hol-system-pi}
q p_i - t \tilde g_i(f_1-\partial_1,f_2-\partial_2,f_3-\partial_3) , \ i=1,2,3 , \\
\label{eq:hol-system-q}
q^2 (\partial_q-f_q) + t \tilde g(f_1-\partial_1,f_2-\partial_2,f_3-\partial_3) , \\
\label{eq:hol-system-t}
q \partial_t - \tilde g(f_1-\partial_1,f_2-\partial_2,f_3-\partial_3) .
\end{gather}
One could prove as in~\cite{BrochetChyzakLairez-2025-FMI}
that the system given by \eqref{eq:hol-system-pi} and~\eqref{eq:hol-system-q} describes
a holonomic module over $\bQ(t)[p_1,p_2,p_3]\langle\partial_1,\partial_2,\partial_3\rangle$
closed by a derivation with respect to~$t$ that is described by~\eqref{eq:hol-system-t}.
But for the purpose of the proof it is sufficient to observe
that the subsequent calculation manifestly terminates outputting meaningful information.

Consider the left $\cD$-ideal generated by the five skew polynomials
in \eqref{eq:hol-system-pi}, \eqref{eq:hol-system-q}, and~\eqref{eq:hol-system-t}.
By construction, any of the two integration algorithms described in~\cite{BrochetChyzakLairez-2025-FMI},
when applied to this left ideal,
will output a nonzero element of~$\bQ(t)\langle\partial_t\rangle$
that annihilates~$\res_{p,q}(F \cL(G)) = \diag_{q,t}(q S(q, t))$,
provided it terminates.
Using Hadrien Brochet's Julia implementation of those algorithms~\cite{Brochet-2025-MCT},
the command~\texttt{MCT} returns after a dozen seconds the skew polynomial
\begin{multline}
\label{eq:lodo-after-diag}
(2624t^{12} - 1752t^{10} - 2160t^8 + 126t^6 + 297t^4) \, \partial_t^5 \\
+ (-15744t^{13} + 21008t^{11} + 10320t^9 - 20340t^7 - 7866t^5 + 1188t^3) \, \partial_t^4 \\
+ (23616t^{14} - 72512t^{12} - 24464t^{10} + 126768t^8 - 432t^6 - 35550t^4 - 2079t^2) \, \partial_t^3 \\
+ (107584t^{13} + 23616t^{11} - 238176t^9 + 68112t^7 + 131004t^5 + 22680t^3 + 2673t) \, \partial_t^2 \\
+ (20992t^{14} + 1728t^{12} - 54336t^{10} - 173952t^8 - 94608t^6 - 44604t^4 - 22680t^2 - 2673) \, \partial_t \\
+ (-10496t^{15} + 29312t^{13} - 36672t^{11} - 142848t^9 - 170064t^7 - 55296t^5) .
\end{multline}
This output is in $t^{-1}\bQ[t^2]\langle t\partial_t \rangle$,
which the expert eye immediately recognizes to correspond to a recurrence relation
between $c_n$, $c_{n+2}$, $c_{n+4}$, etc, for $c_n = s_{n-1,n}/n!$,
thus resembling the conjectured relation~\eqref{eq:conj-rec}.
Upon effectively converting the implied ODE,
one gets a linear relation between $c_n$, $c_{n+2}$, \dots, $c_{n+16}$,
which one next converts to the following linear recurrence relation for $a_k = (2k)! \, c_{2k}$,
satisfied by construction for~$k\geq0$:
\begin{multline}\label{eq:new-rec}
- 83968(2k+15)(2k+13)(2k+11)(9+2k)(2k+7)(2k+5)(2k+3)(2k+1) \\
  \shoveright{\times (k+5)(k+4)(k+3)(k+3)(k+1) \, a_k} \\
\shoveleft{\qquad + 512(2k+15)(2k+13)(2k+11)(2k+9)(2k+7)(2k+5)(2k+3)} \\
  \shoveright{\times (k+5)(k+4)(k+3)(k+2)(328k+557) \, a_{k+1}} \\
\shoveleft{\qquad + 128(2k+15)(2k+13)(2k+11)(9+2k)(7+2k)(5+2k)(5+k)(4+k)(3+k)} \\
  \shoveright{\times (2952k^3 + 20008k^2 + 42776k + 28563) \, a_{k+2}} \\
\shoveleft{\qquad - 512(2k+11)(9+2k)(7+2k)(5+k)(4+k)(2k+15)(2k+13)} \\
  \shoveright{\times (492k^4 + 5561k^3 + 22950k^2 + 41191k + 27597) \, a_{k+3}} \\
\shoveleft{\qquad + 8(2k+15)(2k+13)(2k+11)(9+2k)(5+k)} \\
  \shoveright{\times (5248k^5 + 99728k^4 + 726632k^3 + 2488192k^2 + 3834156k + 1862955) \, a_{k+4}} \\
\shoveleft{\qquad - 48(2k+15)(2k+13)(2k+11)(6+k)} \\
  \shoveright{\times (292k^4 + 3228k^3 + 7085k^2 - 25401k - 82842) \, a_{k+5}} \\
\shoveleft{\qquad - 72(k+7)(2k+15)(2k+13)(120k^3 + 2005k^2 + 10616k + 17805) \, a_{k+6}} \\
+ 18(2k+15)(k+8)(14k^2 - 311k - 2403) \, a_{k+7} + 297(k+9)(k+7) \, a_{k+8} = 0 .
\end{multline}

Finally, note the positivity for all~$k\geq0$ of the leading coefficient of~\eqref{eq:conj-rec},
so that any~$a_0,\dots,a_4$ uniquely determine a solution.
In particular, $a_{k+5},\dots,a_{k+8}$ can be expressed as linear expressions in terms of $a_k,\dots,a_{k+4}$,
and these expressions are well-defined for all~$k\geq0$.
Substituting into~\eqref{eq:new-rec} reveals that any solution to~\eqref{eq:conj-rec} is also a solution to~\eqref{eq:new-rec}.
Consider Sequence A339987, which we have proven to be a solution to~\eqref{eq:new-rec}.
It starts with the values provided by~\eqref{eq:first-values},
which are easily shown to be compatible with~\eqref{eq:conj-rec} and therefore determine a unique solution of~\eqref{eq:conj-rec}.
The leading coefficient of~\eqref{eq:new-rec} is also positive for all~$k\geq0$,
so the same initial conditions determine a unique solution of~\eqref{eq:new-rec}.
The two solutions of \eqref{eq:conj-rec} and~\eqref{eq:new-rec} therefore match
and we have shown that Sequence A339987 satisfies~\eqref{eq:conj-rec} for~$k\geq0$.

A script supporting the complete proof by residues is available in directory \texttt{residues/} of the online appendix.

\section{Proof by graphic recurrences}
\label{sec:graphs}

In the late 1960s, Read developed a method based on combinatorial recurrence
for counting 3-regular graphs~\cite{Read-1970-SUE},
which was later simplified and extended by Wormald~\cite{Wormald-1979-ELG}.
In this section, which is very much inspired by Wormald's work,
we adapt the approach to count the kind of graphs we are interested in,
with only allowed degrees 3 and~1.

In order to do so,
we need to introduce more general graphs whose degrees are in~$\{1,2,3\}$.
More precisely, given non-negative integers $a$, $b$, and~$c$,
we introduce the class~$\cG_{a,b,c}$ of (possibly non-connected) graphs
with exactly $a$~vertices of degree~1, $b$~vertices of degree~2, $c$~vertices of degree~3,
and, for a reason that will become clear,
without any connected component consisting of just two unary vertices linked by one edge.
In what follows we restrict~$b$ to be between~0 and~2.
Additionally, we consider the subclass~$\tilde\cG_{a,2,c}$
of those elements of~$\cG_{a,2,c}$
whose two binary vertices are not adjacent.
Let also $\cU$ denote the class consisting of the single graph
consisting of two unary vertices connected by an edge and labelled 1 and~2.
To address the conjecture,
it proves useful to introduce the classes
\begin{equation}\label{eq:A-from-G0}
\cA_k = \bigcup_{\ell=0}^{\lfloor(k+1)/2\rfloor} \cG_{k+1-2\ell,0,k-1} \times \text{\sc Set}_\ell(\cU) ,
\end{equation}
where disjoint union, product, and {\sc Set} construction are the usual
operations of analytic combinatorics
\cite{FlajoletSedgewick-2009-AC}
(see alternatively the theory of species, e.g.~\cite{BergeronLabelleLeroux-1998-CST,BergeronLabelleLeroux-2013-ITS}).
In particular, $\text{\sc Set}_\ell$~constructs sets of cardinality~$\ell$.
The class~$\cA_k$ contains only graphs that have $2k-1$~edges and $2k$~vertices.

The classes $\cG_{a,0,c}$, $\cG_{a,1,c}$, $\cG_{a,2,c}$, and~$\tilde\cG_{a,2,c}$,
with $a$ and~$c$ ranging over~$\bZ_{\geq0}$,
are interrelated.
We develop a set of reductions that provide a combinatorial recursion between them.
In what follows, some of the unary; binary; cubic vertices
are named $\alpha$, $\alpha'$, etc.; $\beta$, $\beta'$, etc.; $\gamma$, $\gamma'$, etc.; respectively.
The combinatorial recursion to be developed involves several elementary transformations:
the \emph{removal} of a vertex from a graph consists of the deletion of the vertex and all its incident edges;
the \emph{suppression} of a binary vertex consists of removing the vertex
before joining the two vertices formerly adjacent to it;
a graph on $n$~vertices is \emph{improperly labelled}
if its labelling is by two-by-two distinct integers,
without the constraint that the collection of labels is exactly $\{1,\dots,n\}$;
\emph{compressing} an improperly labelled graph on $n$~vertices,
or \emph{compressing} an improper labelling when the graph is clear from the context, consists
in relabelling the vertices in the same total order so as to use exactly all labels from $\{1,\dots,n\}$;
the \emph{unlabelled structure} of a labelled graph
is the unlabelled graph obtained by forgetting the labels;
given a labelled graph~$G$ and an unlabelled graph~$H$ such that
the unlabelled structure of~$G$ is a subgraph of~$H$,
\emph{prolonging} the labelled graph~$G$ with respect to the unlabelled graph~$H$,
often clear from the context,
is the operation of labelling~$H$ in all possible ways that produce a labelled graph~$H'$
such that compressing the improper labelling induced by~$H'$ on the unlabelled structure of~$G$
results in the labelled graph~$G$.

Remember that a graph with $a$~unary vertices, $b$~binary vertices, and $c$~cubic vertices has
$\tfrac12 (a+2b+3c)$~edges.
Introduce the numbers $g_{a,0,c}$, $g_{a,1,c}$, $g_{a,2,c}$, and~$\tilde g_{a,2,c}$,
with $a$ and~$c$ ranging over~$\bZ_{\geq0}$,
that count the classes of graphs $\cG_{a,0,c}$, $\cG_{a,1,c}$, $\cG_{a,2,c}$, and~$\tilde\cG_{a,2,c}$,
respectively,
and likewise introduce the exponential generating functions
\begin{align}
\label{eq:def-G_0}
G_0(q,t) &= \sum_{a,c\geq0} g_{a,0,c} q^{\frac{a+3c}{2}} \frac{t^{a+c}}{(a+c)!} , \\
\label{eq:def-G_1}
G_1(q,t) &= \sum_{a,c\geq0} g_{a,1,c} q^{\frac{a+2+3c}{2}} \frac{t^{a+c+1}}{(a+c+1)!} , \\
\label{eq:def-G_2}
G_2(q,t) &= \sum_{a,c\geq0} g_{a,2,c} q^{\frac{a+4+3c}{2}} \frac{t^{a+c+2}}{(a+c+2)!} , \\
\label{eq:def-tG_2}
\tilde G_2(q,t) &= \sum_{a,c\geq0} \tilde g_{a,2,c} q^{\frac{a+4+3c}{2}} \frac{t^{a+c+2}}{(a+c+2)!} .
\end{align}
\noeqref{eq:def-G_0,eq:def-tG_2}%

To ease the following presentation,
we extend the sequences we just introduced
so that $g_{a,0,c}$, $g_{a,1,c}$, $g_{a,2,c}$, and~$\tilde g_{a,2,c}$
are all zero if any of $a$ and~$c$ is negative.
We also implicitly replace the factors $1/(a+c+i)!$ for~$0\leq i\leq2$ in the generating functions
with piecewise-defined functions that are zero if the argument~$a+c+i$ is negative:
so, we implicitly use the “inverse factorial”~\eqref{eq:inverse-factorial}
that we used in Section~\ref{sec:triplesum}.
But to make the present section more self-contained and written in a more traditional way,
we refrain from using that inverse factorial
and we stick to the common abuse of notation.
The interest of this extension of the sequences is
that the sums defining the generating functions \eqref{eq:def-G_0}--\eqref{eq:def-tG_2}
can be understood as summations over all $a$ and~$c$ in~$\bZ$,
and that the recurrences that we will derive hold without any restrictions on the integers $a$ and~$c$.

The proof is by cases and requires a study of 19~subcases, arranged in a few levels of nesting.
To ease the reading, we have numbered mutually exclusive subcases in subsequent order,
but in doing so we have disregarded the different levels of nesting.

Besides a plain graph, the proof often requires to consider
a graph with a distinguished edge or vertex.
We will speak of \emph{pointing} an edge or vertex,
and of a \emph{pointed} edge or vertex.
According to our notation, a graph of~$\cG_{a,b,c}$
corresponds to $a+b+c$ variants with a pointed vertex
and to $\tfrac12 (a+2b+3c)$ variants with a pointed edge.
This operation of pointing is what introduces derivatives in the equations.

All cases but Case~0' are generalizations from Wormald's presentation,
although the presence of unary vertices introduces many more subcases.
In particular, Case~0 proceeds by pointing and removing an edge:
the counterpart in~\cite{Wormald-1979-ELG} obtains a differential equation with respect to the variable counting vertices
because in cubic graphs, the numbers of edges and vertices are proportional,
but in the present work, this remains a differential equation with respect to~$q$,
which marks edges.
For its part, Case~0' has no counterpart in~\cite{Wormald-1979-ELG}:
it is needed to have derivatives with respect to~$t$
and proceeds by pointing a vertex.
But to avoid many and difficult subcases,
it proved useful to point only unary vertices.

The proof of all the cases starts by fixing a graph from some class, like~$\cG_{a,0,c}$ for Case~0.
The reader will verify that if the class is empty
(e.g., if $a<0$ or~$c<0$, but other situations may force the class to be empty),
implying that no graph can be picked and that no branch of the reasoning by cases (e.g., 0a.\ and~0b.) can be used,
the resulting recurrence relation (e.g.,~\eqref{eq:first-rec-for-g_a,0,c}) holds as a tautology of the form~$0=0$.

\paragraph{Case~0: Graphs with no binary vertex, decomposition from a pointed edge.}

Fix a graph of~$\cG_{a,0,c}$ and identify one of its edges as pointed.
This edge cannot joint two unary vertices,
so only two cases are possible:
\begin{itemize}

\item[0a.] If the pointed edge joins a unary vertex~$\alpha$ and a cubic vertex~$\gamma$,
removing~$\alpha$ turns~$\gamma$ into a binary vertex.
Compressing the thus obtained improperly labelled graph results in a graph of~$\cG_{a-1,1,c-1}$.
The labelled graph so obtained can be obtained from $a+c$ different labellings of the same unlabelled structure.
Conversely, a graph of~$\cG_{a-1,1,c-1}$ has a single binary vertex,~$\beta$,
so attaching a new unary vertex to~$\beta$ makes it cubic,
and the corresponding prolongation results in $a+c$ properly labelled graphs.

\item[0b.] If the pointed edge joins two cubic vertices $\gamma$ and~$\gamma'$,
deleting the pointed edge makes $\gamma$ and~$\gamma'$ binary
and no longer adjacent, as multiple edges are forbidden in~$\cG_{a,0,c}$.
This results in a graph of~$\tilde\cG_{a,2,c-2}$.
Conversely, a graph of~$\tilde\cG_{a,2,c-2}$ contains a uniquely defined unordered pair~$\{\beta,\beta'\}$
of binary vertices.
Joining them by an edge and pointing the latter results bijectively in a graph of~$\cG_{a,0,c}$ with one pointed edge.
\end{itemize}

Collecting these two cases leads to the recurrence relation
\begin{equation}\label{eq:first-rec-for-g_a,0,c}
\tfrac12 (a+3c) g_{a,0,c} = (a+c) g_{a-1,1,c-1} + \tilde g_{a,2,c-2} ,
\end{equation}
which, by the discussion before Case~0, holds for all $a$ and~$c$ in~$\bZ$,
not just for the non-negative $a$ and~$c$ that make~$g_{a,0,c}$ non-zero.
Upon multiplication by $q^{\tfrac12(a+3c)}t^{a+c}/(a+c)!$, we get
\begin{multline*}
\tfrac12 (a+3c) g_{a,0,c} q^{\tfrac12(a+3c)}\frac{t^{a+c}}{(a+c)!} = \\
q (a+c) g_{a-1,1,c-1} q^{\tfrac12(a-2+3c)}\frac{t^{a+c}}{(a+c)!} + q \tilde g_{a,2,c-2} q^{\tfrac12(a-2+3c)}\frac{t^{a+c}}{(a+c)!} .
\end{multline*}
Summing over $a$ and~$c$ in~$\bZ$ now yields
\begin{equation*}
q \partial_q \cdot G_0(q,t) = q t G_1(q,t) + q \tilde G_2(q,t) .
\end{equation*}

\paragraph{Case~0': Graphs with no binary vertex, decomposition from a pointed unary vertex.}

Fix a graph of~$\cG_{a,0,c}$ and identify one of its unary vertices,~$\alpha$, as pointed.
Its neighbor~$\gamma$ must be cubic,
for otherwise the graph would contain a connected component made of two adjacent unary vertices,
which we have excluded.
Removing~$\alpha$ turns~$\gamma$ into a binary vertex.
Compressing the thus obtained improperly labelled graph results in a graph of~$\cG_{a-1,1,c-1}$.
The labelled graph so obtained can be obtained from $a+c$ different labellings of the same unlabelled structure.
Conversely, a graph of~$\cG_{a-1,1,c-1}$ has a single binary vertex,~$\beta$,
so attaching a new unary vertex~$\alpha$ to~$\beta$ makes it cubic,
and after pointing~$\alpha$
the corresponding prolongation results in $a+c$ properly labelled graphs with a pointed unary vertex.

From this simple decomposition, we get the recurrence relation
\begin{equation*}
a g_{a,0,c} = (a+c) g_{a-1,1,c-1} ,
\end{equation*}
valid for all $a$ and~$c$ in~$\bZ$.
Upon multiplication by $q^{\tfrac12(a+3c)}t^{a+c}/(a+c)!$, we get
\begin{equation*}
a g_{a,0,c} \frac{q^{\tfrac12(a+3c)}t^{a+c}}{(a+c)!} = qt g_{a-1,1,c-1} \frac{q^{\tfrac12(a-2+3c)}t^{a+c-1}}{(a+c-1)!} .
\end{equation*}
Observing $\tfrac32 (a+c) - \tfrac12(a+3c) = a$ and
summing over $a$ and~$c$ in~$\bZ$ now yields
\begin{equation*}
\bigl(\tfrac32 t \partial_t - q \partial_q\bigr) \cdot G_0(q,t) = q t G_1(q,t) .
\end{equation*}

\paragraph{Case~1: Graphs with exactly one binary vertex.}

Fix a graph of~$\cG_{a,1,c}$ with unique binary vertex~$\beta$.
Exactly one of the three following situations can happen:
\begin{itemize}

\item[1a.] If $\beta$~is adjacent to two unary vertices,
  thus forming a connected component,
  then the remainder of the graph, after compression,
  is an element of~$\cG_{a-2,0,c}$.
  Conversely, enriching an element of~$\cG_{a-2,0,c}$ with a connected component
  made of a binary vertex and an unordered pair of unary vertices,
  and the corresponding prolongation results in $\tfrac12(a+c+1)(a+c)(a+c-1)$~elements of~$\cG_{a,1,c}$.

\item[1b.] If $\beta$~is adjacent to one unary vertex~$\alpha$ and one cubic vertex~$\gamma$,
  then removing $\alpha$ and~$\beta$ makes~$\gamma$ binary.
  After compressing the thus obtained improperly labelled graph
  results in an element of~$\cG_{a-1,1,c-1}$.
  Conversely, given an element of~$\cG_{a-1,1,c-1}$,
  call~$\beta$ its binary vertex.
  Consider the corresponding unlabelled structure
  and graft a vertex~$\beta'$ to~$\beta$, then a vertex~$\alpha$ to~$\beta'$,
  thus turning~$\beta$ into a cubic vertex.
  The corresponding prolongation results in $(a+c+1)(a+c)$~elements of~$\cG_{a,1,c}$.

\item[1c.] If $\beta$~is adjacent to two cubic vertices $\gamma$ and~$\gamma'$,
  then removing~$\beta$, thus making $\gamma$ and~$\gamma'$ binary,
  and finally compressing the thus obtained improperly labelled graph
  results in a graph of~$\cG_{a,2,c-2}$.
  The same resulting graph can be obtained from $a+c+1$~choices in~$\cG_{a,1,c}$.

\end{itemize}
Collecting contributions leads to the recurrence relation
\begin{multline*}
g_{a,1,c} = \tfrac12(a+c+1)(a+c)(a+c-1) g_{a-2,0,c} \\
  {} + (a+c+1)(a+c) g_{a-1,1,c-1} + (a+c+1) g_{a,2,c-2} ,
\end{multline*}
valid for all $a$ and~$c$ in~$\bZ$.
Upon multiplication by $q^{\tfrac12(a+2+3c)}t^{a+c+1}/(a+c+1)!$, we get
\begin{multline*}
g_{a,1,c} q^{\tfrac12(a+2+3c)}\frac{t^{a+c+1}}{(a+c+1)!} =
  \tfrac12 q^2 t^3 g_{a-2,0,c} q^{\tfrac12(a-2+3c)}\frac{t^{a+c-1}}{(a+c-1)!} \\
  {} + q^2 t^2 g_{a-1,1,c-1} q^{\tfrac12(a-2+3c)}\frac{t^{a+c-1}}{(a+c-1)!}
  {} + q^2 t g_{a,2,c-2} q^{\tfrac12(a-2+3c)}\frac{t^{a+c}}{(a+c)!} .
\end{multline*}
Summing over $a$ and~$c$ in~$\bZ$ now yields
\begin{equation*}
G_1(q,t) = \tfrac12 q^2 t^3 G_0(q,t) + q^2 t^2 G_1(q,t) + q^2 t G_2(q,t) .
\end{equation*}

\paragraph{Case~2: Graphs with exactly two binary vertices.}

Fix a graph of~$\cG_{a,2,c}$ with two binary vertices, $\beta$ and~$\beta'$.
A first option is that $\beta$ and~$\beta'$ are not adjacent, in which case:
\begin{itemize}

\item[2a.] The graph is just an element of the subclass~$\tilde\cG_{a,2,c}$.

\end{itemize}
Otherwise, $\beta$~is also connected to a vertex~$\nu$ distinct from~$\beta'$
and $\beta'$~is also connected to a vertex~$\nu'$ distinct from~$\beta$.
If $\nu$ and~$\nu'$ are distinct, three exclusive cases are possible:
\begin{itemize}

\item[2b.] $\nu$ and~$\nu'$ are unary vertices,
  so that $\nu$, $\beta$, $\beta'$, and~$\nu'$ form a connected component.
  The remainder of the graph, after compression, is just an element of~$\cG_{a-2,0,c}$.
  Conversely, an element of~$\cG_{a-2,0,c}$ can be prolonged into an element of~$\cG_{a,2,c}$
  in $\binom{a+c+2}{2} (a+c)(a+c-1)$ ways
  by choosing an unordered pair of labels for the two added binary vertices,
  then choosing a label for the unary vertex added to the binary vertex with smaller label,
  and finally choosing a label for the unary vertex added to the binary vertex with larger label.

\item[2c.] $\nu$~is a unary vertex and $\nu'$~is a cubic vertex.
  Then, removing $\nu$, $\beta$, and~$\beta'$ results, after compression,
  in an element of~$\cG_{a-1,1,c-1}$.
  Conversely, given an element of~$\cG_{a-1,1,c-1}$,
  an element of~$\cG_{a,2,c}$ is obtained by attaching
  a binary vertex~$\beta'$ to the only binary vertex of the original graph,
  then a binary vertex~$\beta$ to~$\beta'$,
  then a unary vertex~$\alpha$ to~$\beta$.
  This prolongation results in $(a+c+2)(a+c+1)(a+c)$ labelled graphs
  as there are that many choices of indices for the ordered list $\beta'$, $\beta$, $\alpha$.

\item[2d.] $\nu$ and~$\nu'$ are cubic vertices,
  so that removing $\beta$ and~$\beta'$ and compressing the thus obtained improperly labelled graph
  results in an element of~$\cG_{a,2,c-2}$.
  Conversely, given an element of~$\cG_{a,2,c-2}$,
  call $\gamma$ and~$\gamma'$ its two binary vertices, with the convention that $\gamma$~has the smaller of the two labels.
  Then, graphs of~$\cG_{a,2,c}$ can be obtained by inserting two binary vertices
  so as to form a chain $\gamma$, $\beta$, $\beta'$, $\gamma'$, adjacent in that order.
  The corresponding prolongation results in $(a+c+2)(a+c+1)$ labelled graphs,
  as there are that many choices of possible labels for $\beta$ and~$\beta'$.

\end{itemize}
If $\nu$ and~$\nu'$ are equal, then $\beta$, $\beta'$, and~$\nu$ form a 3-cycle.
Because $\nu$~has degree at least~$2$ and is neither $\beta$ nor~$\beta'$, it must be cubic.
Therefore, $\nu$~is connected to a third, distinct vertex,~$\nu''$.
Only two exclusive cases are thus possible:
\begin{itemize}

\item[2e.] The vertex~$\nu''$ is unary,
  so that $\nu$, $\beta$, $\beta'$, and~$\nu''$ form a connected component.
  The remainder of the graph, after compression, is just an element of~$\cG_{a-1,0,c-1}$.
  Conversely, an element of~$\cG_{a-1,0,c-1}$ can be prolonged into an element of~$\cG_{a,2,c}$
  in $\binom{a+c+2}{2} (a+c) (a+c-1)$ ways
  by choosing an unordered pair of labels for the two added binary vertices,
  a label for the added cubic vertex, and a label for the unary vertex,
  before relabelling the original graph.

\item[2f.] The vertex~$\nu''$ is cubic.
  Removing $\beta$, $\beta'$, and~$\nu$, thus making~$\nu''$ binary,
  then compressing the thus obtained improperly labelled graph
  leads to an element of~$\cG_{a,1,c-2}$.
  Conversely, given an element of~$\cG_{a,1,c-2}$, graphs of~$\cG_{a,2,c}$ can be obtained
  by attaching a cubic vertex to the binary vertex of the graph,
  then an unordered pair of binary vertices to the new cubic vertex.
  The corresponding prolongation results in $(a+c+2) \binom{a+c+1}2$ ways.

\end{itemize}
Collecting the contributions of all cases, we get the recurrence relation
\begin{multline*}
g_{a,2,c} = \tilde g_{a,2,c} \\
  {} + \tfrac12 (a+c+2)(a+c+1)(a+c)(a+c-1) g_{a-2,0,c} \\
  {} + (a+c+2)(a+c+1)(a+c) g_{a-1,1,c-1} \\
  {} + (a+c+2)(a+c+1) g_{a,2,c-2} \\
  {} + \tfrac12 (a+c+2)(a+c+1)(a+c)(a+c-1) g_{a-1,0,c-1} \\
  {} + \tfrac12 (a+c+2)(a+c+1)(a+c) g_{a,1,c-2} ,
\end{multline*}
valid for all $a$ and~$c$ in~$\bZ$.
Upon multiplication by $q^{\tfrac12(a+4+3c)}t^{a+c+2}/(a+c+2)!$, we get
\begin{multline*}
g_{a,2,c} q^{\tfrac12(a+4+3c)}\frac{t^{a+c+2}}{(a+c+2)!} = \tilde g_{a,2,c} q^{\tfrac12(a+4+3c)}\frac{t^{a+c+2}}{(a+c+2)!} \\
  {} + \tfrac12 q^3 t^4 g_{a-2,0,c} q^{\tfrac12(a-2+3c)}\frac{t^{a+c-2}}{(a+c-2)!} \\
  {} + q^3 t^3 g_{a-1,1,c-1} q^{\tfrac12(a-2+3c)}\frac{t^{a+c-1}}{(a+c-1)!} \\
  {} + q^3 t^2 g_{a,2,c-2} q^{\tfrac12(a-2+3c)}\frac{t^{a+c}}{(a+c)!} \\
  {} + \tfrac12 q^4 t^4 g_{a-1,0,c-1} q^{\tfrac12(a-4+3c)}\frac{t^{a+c-2}}{(a+c-2)!} \\
  {} + \tfrac12 q^4 t^3 g_{a,1,c-2} q^{\tfrac12(a-4+3c)}\frac{t^{a+c-1}}{(a+c-1)!} .
\end{multline*}
Summing over $a$ and~$c$ in~$\bZ$ now yields
\begin{multline*}
G_2(q,t) = \tilde G_2(q,t) + \tfrac12 q^3 t^4 G_0(q,t) + q^3 t^3 G_1(q,t) + q^3 t^2 G_2(q,t) \\
  {} + \tfrac12 q^4 t^4 G_0(q,t) + \tfrac12 q^4 t^3 G_1(q,t) .
\end{multline*}

\paragraph{Case~3: Graphs with exactly two binary vertices, assumed to not be adjacent.}

Focusing on the class of graphs with no non-adjacent binary vertices
makes us consider a subclass of the preceding class,
but we have to address it for itself as it is a sub-product of two of the earlier transformations.
Fix a graph of~$\tilde\cG_{a,2,c}$ with two binary vertices, one of which is pointed.
Call~$\beta$ this pointed vertex and $\beta'$ the other binary vertex.
The vertex~$\beta$ is connected to two other vertices, $\nu$ and~$\nu'$.

If $\nu$ and~$\nu'$ are not adjacent, two exclusive cases have to be considered:
\begin{itemize}

\item[3a.] If both $\nu$ and~$\nu'$ are unary, then $\nu$, $\beta$, and~$\nu'$ form a connected component.
  Removing the whole connected component and compressing the labelling of the obtained graph
  results in an element of~$\cG_{a-2,1,c}$.
  Conversely, given an element of~$\cG_{a-2,1,c}$,
  its unlabelled structure can be augmented with an additional unlabelled connected component
  consisting of two unary vertices attached to the same binary vertex.
  The corresponding prolongation into graphs of~$\tilde\cG_{a,2,c}$ with a pointed binary vertex
  results in $(a+c+2) \binom{a+c+1}2$ elements, as one has to choose
  a label for the pointed binary vertex and an unordered pair of labels for the unary vertices.

\item[3b.] If $\nu$ and~$\nu'$ are not both unary,
  then suppressing~$\beta$ and compressing the thus obtained labelling results
  in an element of~$\cG_{a,1,c}$ with a pointed edge other than the two edges incident to~$\beta'$.
  Conversely, consider an element of~$\cG_{a,1,c}$ with a pointed edge that is not incident to the binary vertex.
  Inserting a second binary vertex into the pointed edge,
  pointing it,
  and prolonging the original graph with respect to the obtained unlabelled structure
  results in $a+c+2$ elements of~$\tilde\cG_{a,2,c}$ with a pointed binary vertex.

\end{itemize}

Otherwise, $\nu$ and~$\nu'$ are adjacent, and therefore cubic
as their degree is at least~2 and cannot be equal to~$\beta'$, which is not adjacent to~$\beta$.
Therefore,
the vertex~$\nu$ is connected to a third vertex~$\mu$, distinct from $\beta$ and~$\nu'$,
and likewise the vertex~$\nu'$ is connected to a third vertex~$\mu'$, distinct from $\beta$ and~$\nu$.

\begin{itemize}

\item[3c.] Assume that $\mu$ and~$\mu'$ are distinct.
Then, deleting the edge between $\nu$ and~$\nu'$, then suppressing $\nu$ and~$\nu'$,
keeping $\beta$ pointed,
and finally compressing the thus obtained labelling
results in an element of~$\cG_{a,2,c-2}$ with a pointed binary vertex.
Conversely, consider an element of~$\cG_{a,2,c-2}$ with a pointed binary vertex~$\beta$.
Inserting cubic vertices into both edges incident to~$\beta$,
joining those two vertices,
keeping $\beta$ pointed,
and considering the corresponding prolongation
results in $(a+c+2)(a+c+1)$ elements of~$\tilde\cG_{a,2,c}$ with a pointed binary vertex.

\end{itemize}

Otherwise, $\mu$ and~$\mu'$ are the same vertex.
There are three exclusive cases for the relation between $\mu$ and~$\beta'$:
\begin{itemize}

\item If $\mu$~is equal to~$\beta'$, then:
  \begin{itemize}
  \item[3d.] $\beta$, $\nu$, $\nu'$, and~$\beta'$ form a connected component.
  Removing the whole connected component and compressing the labelling of the obtained graph
  results in an element of~$\cG_{a,0,c-2}$.
  Conversely, given an element of~$\cG_{a,0,c-2}$,
  its unlabelled structure can be augmented with an additional unlabelled connected component
  consisting of two binary vertices, both attached to the same two cubic vertices,
  the latter being connected to realize their degree~3.
  The corresponding prolongation into graphs of~$\tilde\cG_{a,2,c}$ with a pointed binary vertex
  results in $(a+c+2)(a+c+1) \binom{a+c}2$ elements, as one has to choose:
  a label for the pointed binary vertex, then a label for the other binary vertex,
  and finally an unordered pair of labels for the cubic vertices.
  \end{itemize}

\item If $\mu$~is connected to a third vertex, distinct from $\nu$ and~$\nu'$,
  which happens to be~$\beta'$,
  then $\beta'$~is connected to a second vertex,~$\mu''$, distinct from~$\mu$.
  Then:
  \begin{itemize}
  \item[3e.] If $\mu''$~is unary, then $\beta$, $\nu$, $\nu'$, $\mu$, $\beta'$, and~$\mu''$ form a connected component.
    Removing the whole connected component and compressing the labelling of the obtained graph
    results in an element of~$\cG_{a-1,0,c-3}$.
    Conversely, given an element of~$\cG_{a-1,0,c-3}$,
    its unlabelled structure can be augmented with an additional unlabelled connected component
    consisting of a unary vertex attached to a binary vertex,
    itself attached to a cubic vertex,
    whose two other adjacent vertices are cubic, adjacent to one another,
    and share a common third adjacent vertex that is binary and pointed.
    The corresponding prolongation into graphs of~$\tilde\cG_{a,2,c}$ with a pointed binary vertex
    results in $(a+c+2) \binom{a+c+1}2 (a+c-1)(a+c-2)(a+c-3)$ elements, as one has to choose:
    a label for the pointed binary vertex,
    next an unordered pair of labels for the cubic vertices attached to this pointed vertex,
    next a label for the third cubic vertex,
    then a label for the other binary vertex,
    and finally a label for the unary vertex.
  \item[3f.] If $\mu''$~is not unary, then it is cubic.
    Then, removing $\beta$, $\nu$, $\nu'$, $\mu$, and~$\beta'$,
    thus making~$\mu''$ binary,
    results in an element of~$\cG_{a,1,c-4}$.
    Conversely, given an element of~$\cG_{a,1,c-4}$,
    its unlabelled structure can be augmented
    by attaching to its only binary vertex a new binary vertex,
    then a cubic vertex to this new vertex,
    then two other cubic vertices adjacent to this first cubic vertex
    in such a way that the two cubic vertices are also adjacent to one another
    and share a common third adjacent vertex that is binary and pointed.
    The corresponding prolongation into graphs of~$\tilde\cG_{a,2,c}$ with a pointed binary vertex
    results in $(a+c+2) \binom{a+c+1}2 (a+c-1)(a+c-2)$ elements, as one has to choose:
    a label for the pointed binary vertex,
    next an unordered pair of labels for the cubic vertices attached to this pointed vertex,
    next a label for the third cubic vertex,
    and finally a label for the other binary vertex.
  \end{itemize}

\item If $\mu$~is connected to a third vertex,~$\mu''$, distinct from $\nu$, $\nu'$, and~$\beta'$:
  \begin{itemize}
  \item[3g.] If $\mu''$~is unary, then $\beta$, $\nu$, $\nu'$, $\mu$, and~$\mu''$ form a connected component.
    Removing the whole connected component and compressing the labelling of the obtained graph
    results in an element of~$\cG_{a-1,1,c-3}$.
    Conversely, given an element of~$\cG_{a-1,1,c-3}$,
    its unlabelled structure can be augmented with an additional unlabelled connected component
    consisting of a unary vertex attached to a cubic vertex,
    whose two other adjacent vertices are cubic, adjacent to one another,
    and share a common third adjacent vertex that is binary and pointed.
    The corresponding prolongation into graphs of~$\tilde\cG_{a,2,c}$ with a pointed binary vertex
    results in $(a+c+2) \binom{a+c+1}2 (a+c-1)(a+c-2)$ elements, as one has to choose:
    a label for the pointed binary vertex,
    next an unordered pair of labels for the cubic vertices attached to this pointed vertex,
    next a label for the third cubic vertex,
    and finally a label for the unary vertex.
  \item[3h.] If $\mu''$~is not unary, then it is cubic.
    Removing $\beta$, $\nu$, $\nu'$, and~$\mu$, thus making~$\mu''$ binary,
    then compressing the thus obtained labelling,
    and finally pointing~$\mu''$
    results in an element of~$\cG_{a,2,c-4}$ with a pointed binary vertex.
    Conversely, given an element of~$\cG_{a,2,c-4}$ with a pointed binary vertex,
    its unlabelled structure can be augmented
    by attaching to its pointed binary vertex a new cubic vertex,
    then two other cubic vertices adjacent to this first cubic vertex
    in such a way that the two cubic vertices are also adjacent to one another
    and share a common third adjacent vertex that is binary and pointed.
    The corresponding prolongation into graphs of~$\tilde\cG_{a,2,c}$ with a pointed binary vertex
    results in $(a+c+2) \binom{a+c+1}2 (a+c-1)$ elements, as one has to choose:
    a label for the pointed binary vertex,
    next an unordered pair of labels for the cubic vertices attached to this pointed vertex,
    and finally a label for the third cubic vertex.
  \end{itemize}

\end{itemize}
Collecting the contributions of all cases, we get the recurrence relation
\begin{multline*}
2 \tilde g_{a,2,c} = \tfrac12 (a+c+2)(a+c+1)(a+c) \times g_{a-2,1,c} \\
  {} + (a+c+2) \times \bigl(\tfrac12 (a+2+3c) -2\bigr) g_{a,1,c} \\
  {} + (a+c+2)(a+c+1) \times 2 g_{a,2,c-2} \\
  {} + \tfrac12 (a+c+2)(a+c+1)(a+c)(a+c-1) \times g_{a,0,c-2} \\
  {} + \tfrac12 (a+c+2)(a+c+1)(a+c)(a+c-1)(a+c-2)(a+c-3) \times g_{a-1,0,c-3} \\
  {} + \tfrac12 (a+c+2)(a+c+1)(a+c)(a+c-1)(a+c-2) \times g_{a,1,c-4} \\
  {} + \tfrac12 (a+c+2)(a+c+1)(a+c)(a+c-1)(a+c-2) \times g_{a-1,1,c-3} \\
  {} + \tfrac12 (a+c+2)(a+c+1)(a+c)(a+c-1) \times 2 g_{a,2,c-4} ,
\end{multline*}
valid for all $a$ and~$c$ in~$\bZ$.
Upon multiplication by $2q^{\tfrac12(a+4+3c)}t^{a+c+2}/(a+c+2)!$, we get
\begin{multline*}
4 \tilde g_{a,2,c} \frac{q^{\tfrac12(a+4+3c)}t^{a+c+2}}{(a+c+2)!} =
  q^2 t^3 g_{a-2,1,c} \frac{q^{\tfrac12(a+3c)}t^{a+c-1}}{(a+c-1)!} \\
  {} + q t (a-2+3c) g_{a,1,c} \frac{q^{\tfrac12(a+2+3c)}t^{a+c+1}}{(a+c+1)!} \\
  {} + 4 q^3 t^2 g_{a,2,c-2} \frac{q^{\tfrac12(a-2+3c)}t^{a+c}}{(a+c)!}
  {} + q^5 t^4 g_{a,0,c-2} \frac{q^{\tfrac12(a-6+3c)}t^{a+c-2}}{(a+c-2)!} \\
  {} + q^7 t^6 g_{a-1,0,c-3} \frac{q^{\tfrac12(a-10+3c)}t^{a+c-4}}{(a+c-4)!}
  {} + q^7 t^5 g_{a,1,c-4} \frac{q^{\tfrac12(a-10+3c)}t^{a+c-3}}{(a+c-3)!} \\
  {} + q^6 t^5 g_{a-1,1,c-3} \frac{q^{\tfrac12(a-8+3c)}t^{a+c-3}}{(a+c-3)!}
  {} + 2 q^6 t^4 g_{a,2,c-4} \frac{q^{\tfrac12(a-8+3c)}t^{a+c-2}}{(a+c-2)!} .
\end{multline*}
Summing over $a$ and~$c$ in~$\bZ$ now yields
\begin{multline}
\label{eq:deqns-tG2-nosimpl}
4 \tilde G_2(q,t) = q^2 t^3 G_1(q,t) + q t (2 q\partial_q - 4) \cdot G_1(q,t) \\
  {} + 4 q^3 t^2 G_2(q,t) + q^5 t^4 G_0(q,t) + q^7 t^6 G_0(q,t) \\
  {} + q^7 t^5 G_1(q,t) + q^6 t^5 G_1(q,t) + 2 q^6 t^4 G_2(q,t) .
\end{multline}

After this analysis, gathering and simplifying the obtained equations
yields a differential system in $G_0(q,t)$, $G_1(q,t)$, $G_2(q,t)$, and~$\tilde G_2(q,t)$:
\begin{align}
\label{eq:deqns-G0}
\partial_q \cdot G_0(q,t) &= t G_1(q,t) + \tilde G_2(q,t) , \\
\label{eq:deqns-G0-bis}
(3 t \partial_t - 2 q \partial_q) \cdot G_0(q,t) &= 2 q t G_1(q,t) , \\
\label{eq:deqns-G1}
G_1(q,t) &= \tfrac12 q^2 t^3 G_0(q,t) + q^2 t^2 G_1(q,t) + q^2 t G_2(q,t) , \\
\label{eq:deqns-G2}
G_2(q,t) &= \tilde G_2(q,t) + \tfrac12 q^3 t^4 G_0(q,t) + q^3 t^3 G_1(q,t) \\
  \notag
  &\qquad + q^3 t^2 G_2(q,t) + \tfrac12 q^4 t^4 G_0(q,t) + \tfrac12 q^4 t^3 G_1(q,t) , \\
\label{eq:deqns-tG2}
4 \tilde G_2(q,t) &= q^2 t^3 G_1(q,t) + 2 q^2 t\partial_q \cdot G_1(q,t) - 4 q t G_1(q,t) \\
  \notag
  &\qquad + 4 q^3 t^2 G_2(q,t) + q^5 t^4 G_0(q,t) + q^7 t^6 G_0(q,t) \\
  \notag
  &\qquad + q^7 t^5 G_1(q,t) + q^6 t^5 G_1(q,t) + 2 q^6 t^4 G_2(q,t) .
\end{align}
\noeqref{eq:deqns-G0-bis,eq:deqns-G1,eq:deqns-G2}

At this point, eliminating $G_1(q,t)$, $G_2(q,t)$, and~$\tilde G_2(q,t)$ between the four equations is easy:
this can be done by a Gröbner-basis calculation for modules
over the Ore algebra~$\cR = \bQ(q,t)\langle\partial_q,\partial_t\rangle$
(in which $\partial_q q = q\partial_q + 1$, $\partial_t t = t\partial_t + 1$,
and all other pairs of generators commute).
To this end,
one encodes \eqref{eq:deqns-G0}--\eqref{eq:deqns-tG2}
as elements of a free left $\cR$-module of rank~4
that we denote $\cR\cdot g_0\oplus \cR\cdot g_1\oplus \cR\cdot g_2\oplus \cR\cdot \tilde g_2$
after choosing suggestive names $g_0, g_1, g_2, \tilde g_2$ for the elements of its canonical basis.
For example, \eqref{eq:deqns-G0}~becomes
$\partial_q\cdot g_0 - t\cdot g_1 - \tilde g_2$,
and likewise for \eqref{eq:deqns-G0-bis}--\eqref{eq:deqns-tG2}.
Choosing a monomial ordering that has $g_1,g_2,\tilde g_2$
lexicographically larger than $g_0,\partial_q,\partial_t$
and that sorts by the total degree in $(\partial_q,\partial_t)$ with~$\partial_q > \partial_t$,
we obtain a Gröbner basis consisting of five generators,
with only two generators in~$\cR\cdot g_0$:
omitting~$g_0$, these are
\begin{align}
\label{eq:elim-G_0-mixed}
&2q(q^6t^4 + 2q^3t^2 - 2q^2t^2 - 2) \partial_q - 3t(q^6t^4 + 2q^3t^2 - 2) \partial_t - 2q^3t^4(q^4t^2 + 1) , \\
\label{eq:elim-G_0-dt-only}
&9t^3q^3(q^6t^4 + 2q^3t^2 - 2q^2t^2 - 2) \partial_t^2 \\
\notag
&\quad + (3q^{15}t^{10} + 18q^{12}t^8 - 24q^{11}t^8 + 18q^{10}t^8 + 9q^9t^6 - 60q^8t^6 + 72q^7t^6 \\
\notag
&\qquad\qquad - 36q^6t^6 - 18q^6t^4 + 18q^5t^4 - 36q^4t^4 - 78q^3t^2 + 24q^2t^2 + 24) \partial_t \\
\notag
&\quad - t^3q^3(q^{15}t^8 - 3q^{13}t^8 + 4q^{12}t^6 - 8q^{11}t^6 - 18q^{10}t^6 + 24q^9t^6 - 9q^8t^6 \\
\notag
&\qquad\qquad - 28q^8t^4 + 18q^7t^4 + 52q^6t^4 - 54q^5t^4 - 8q^6t^2 + 18q^4t^4 \\
\notag
&\qquad\qquad + 40q^5t^2 + 36q^4t^2 - 64q^3t^2 + 34q^2t^2 + 4q^3 + 16) ,
\end{align}
corresponding to two PDE satisfied by~$G_0(q,t)$.
We have not attempted to do the calculation by hand;
a script for this elimination is available in directory \texttt{graphs/} of the online appendix.

One finally has to remember that the classes~$\cG_{a,0,c}$
exclude graphs with connected components in the singleton class~$\cU$.
The connection between the classes~$\cA_k$ we are interested in
and those classes is done by~\eqref{eq:A-from-G0}.
Because the exponential generating function of~$\cU$ is~$qt^2/2$
and because the {\sc Set} construction is reflected by an exponential function,
we have to put back copies of~$\cU$ by considering
$S(q,t) = G_0(q,t) \exp(qt^2/2)$.
A differential system for it is obtained
by right multiplying \eqref{eq:elim-G_0-mixed}--\eqref{eq:elim-G_0-dt-only} with~$\exp(-qt^2/2)$,
making the exponential pass to the left,
then left multiplying with~$\exp(qt^2/2)$,
and doing so we obtain the same system
as the system~\eqref{eq:bivariate-sys-only-dt}--\eqref{eq:bivariate-sys-mixed}
obtained by the generalization of~\cite{ChyzakMishna-2025-DES} suggested in Section~\ref{sec:residues}.
Correspondingly, a diagonal computation to represent~$\diag_{q,t}(q S(q, t))$
leads again to the linear differential operator~\eqref{eq:lodo-after-diag},
from which we prove the conjectured recurrence relation~\eqref{eq:conj-rec} again.

As a remark, let us briefly comment on what happens when we analyze
another transformation with a combinatorial meaning:
we tried reducing graphs with two adjacent binary vertices,
that is, elements of~$\cG_{a,2,c} \setminus \tilde\cG_{a,2,c}$.
Indeed, by suppressing the one of the two binary vertices whose non-binary neighbor is of smaller label,
we obtain a graph of~$\cG_{a,1,c}$, unless the two non-binary neighbors are in fact the same vertex,
leading either to the removal of a 4-vertex connected component if the shared neighbor is unary,
or to another reduction to a graph of~$\cG_{a,1,c-2}$ if it is cubic.
The reader will check the relation
\begin{multline*}
g_{a,2,c} - \tilde g_{a,2,c} = (a+c+2) g_{a,1,c} \\
  {} + \tbinom{a+c+2}2 (a+c) (a+c-1) g_{a-1,0,c-1} + \tbinom{a+c+2}2 (a+c) g_{a,1,c-2} ,
\end{multline*}
and that it leads to
\begin{equation*}
G_2(q,t) - \tilde G_2(q,t) = q t G_1(q,t) + \tfrac12 q^4 t^4 G_0(q,t) + \tfrac12 q^4 t^3 G_1(q,t) .
\end{equation*}
This nice simple equation is however also obtained by removing $qt$~times~\eqref{eq:deqns-G1} from~\eqref{eq:deqns-G2}.

A script supporting the calculations on the differential system \eqref{eq:deqns-G0}--\eqref{eq:deqns-tG2},
first elimination, next twist by~$\exp(qt^2/2)$, and finally comparison to~\eqref{eq:bivariate-sys-only-dt}--\eqref{eq:bivariate-sys-mixed},
is available in directory \texttt{graphs/} of the online appendix.

\section{Conclusion: computer proofs vs. human proofs, really?}
\label{sec:conclusion}

In this concluding section,
we want to comment on the role of computation in our proofs,
and in particular on its relation to the ease of the proof search.

The extensive efforts we have needed
to obtain a complete proof by the triple summation of Section~\ref{sec:triplesum}
may come as a surprise to our readers:
this is paradoxical as creative telescoping has been presented since its introduction as a proof technique to derive recurrences satisfied by parametrized sums.
However, although implementations of creative-telescoping algorithms almost always return correct recurrences,
the calculation they perform is never a proof by itself
and the additional proof steps that a human needs to perform after it
are very prone to logical gaps.
The fact that creative telescoping is used in a context
where the obtained recurrences can be checked numerically explains why such gaps can easily remain unnoticed.
Yet, difficult situations have been observed for at least a decade,
and exemplified several times with multiple sums:
an attempt to combine creative telescoping with computerized formal proofs,
in the sense of proofs mechanically checked by a computer in a proof assistant,
revealed how difficult it is to manipulate symbolically the certificates output by the algorithms,
and that rational-function normalization is insufficient~\cite{ChyzakMahboubiSibutPinoteTassi-2014-CAB};
Koutschan and Wong~\cite{KoutschanWong-2021-CTM} studied some triple binomial sum
to highlight a host of issues that can be encountered and
that the available counter-measures are typically not systematic.
The issues they encountered include:
summations not over natural boundaries as auxiliary sums;
bounds that depend on parameters;
singularities of the certificates.
Our triple sum~\eqref{eq:triple-sum-with-natural-bounds} is just another example of these problems,
possibly made more difficult than it is for the sum in~\cite{KoutschanWong-2021-CTM}
by the presence of a polynomial in the denominator of the summand.
Note that similar problems have in fact also been observed for indefinite sums:
a notable example, $\sum_k \binom{2k-3}{k}/4^k$, is given as \cite[Example~1]{AbramovPetkovsek-2008-OBS}.

After this and with regard to formal proofs,
we could say that conventional creative telescoping alone,
as described in Section~\ref{sec:shortcomings},
cannot be called a general proving method in its form dealing with recurrence relations:
creative-telescoping algorithms do not provide with more than a second guess,
corroborating the original guess based on Hermite--Padé approximants
by Kauers and Koutschan~\cite{KauersKoutschan-2023-SDF}.
This contrasts with the use of the same approach in a differential setting in Section~\ref{sec:residues},
where creative telescoping easily leads to complete formal proofs.
But supplementing the initial guess of a recurrence in Section~\ref{sec:shortcomings} with the a~posteriori validation of Section~\ref{sec:triplesum-ok}
is another incarnation of Pólya's “First guess, then prove”.

Our proof by graphic recurrences presented in Section~\ref{sec:graphs}
resorts to the computer only at its last step:
for an elimination between \eqref{eq:deqns-G0}--\eqref{eq:deqns-tG2} leading to \eqref{eq:elim-G_0-mixed}--\eqref{eq:elim-G_0-dt-only}.
This last calculation could in principle be performed algorithmically by a human.
We would however not want to be like a fox who effaced its tracks:
although the final proof is complete with little effect of a computer,
we did use a computer for the proof search.
Being aware of~\cite{Wormald-1979-ELG},
it was expected that the combinatorial approach should adapt to $(3,1)$-regular graphs,
but also obvious that this would demand technical additions.
As a matter of fact, our first totally human attempts to write a proof resulted
in wrong equations that we did not fix by pure thinking.
Instead, we devised a reversed approach:
we implemented a known backtracking algorithm~\cite[Section~2]{ReadWormald-1981-C10}
to count graphs with a given degree sequence
and used it to generate truncations of $G_0$, $G_1$, $G_2$ up to 200~vertices
(we did not try to take the non-adjacency of degree-2 vertices in the generation, which would have been needed for~$\tilde G_2$);
we could next use linear-algebra-based guessing
to fix the polynomial coefficients appearing in \eqref{eq:deqns-G0}--\eqref{eq:deqns-tG2}, thus identifying computation mistakes and reasoning errors;
the proof in Section~\ref{sec:graphs} could finally be adjusted
into its present form.

It is interesting that the monomials~$q^\alpha t^\beta$ that appear in the system
have a strong combinatorial meaning:
they indicate the shape of the subgraph that is removed in the transformation,
and often enough suggest a connected component to be removed.
The fixes made possible by the calculation are the following:
\begin{itemize}
\item Calculation/copy-paste mistakes led to an extra $\partial_t$ beside a~$t$ (case~0a.) and in some~$\tfrac14$ instead of a~$\tfrac12$ (case~2e.).
\item A confusion between ordered and unordered pairs (case~2b.) led to another~$\tfrac14$ instead of~$\tfrac12$.
\item In the case of graphs with exactly two binary vertices, assumed to not be adjacent,
  for the subcase when the vertices $\nu$ and~$\nu'$ are not both unary (case~3b.),
  we had overlooked that the suppression of the binary vertex~$\beta$
  could not make the resulting edge adjacent to~$\beta'$,
  hence the term $2q\partial_q - 4 = 2(q\partial_q - 2)$, not our original~$2q\partial_q$,
  in~\eqref{eq:deqns-tG2-nosimpl};
  as $q\partial_q$~is used to point edges,
  a term~$q\partial_q - 2$ meant pointing all but two edges,
  from which the intuition of the missing argument was immediate.
\item Again in the case of graphs with exactly two binary vertices, assumed to not be adjacent,
  we had overlooked the special case when the vertices $\nu$ and~$\nu'$ are both unary (case~3a.),
  thus making $\nu$, $\beta$, and~$\nu'$ form a connected component
  in which the suppression of~$\beta$ would lead to a forbidden structure:
  we had erroneously not distinguished the two situations “both unary” and “not both unary”,
  which was revealed by a missing term~$q^2 t^3 G_1(q,t)$;
  intuiting that this term could represent a Cartesian product
  led to analyze what graph class has~$q^2 t^3$ as a generating function,
  immediately spotting the overlooked possibility of a connected component from the forbidden class~$\cU$.
\end{itemize}

The readers might also find interesting to know
that researching the proof by residues (Section~\ref{sec:residues}) was carried out in less than two days,
while the approach by graphic recurrences (Section~\ref{sec:graphs}) needed weeks of full-time work,
and the approach by a triple sum (Section~\ref{sec:triplesum}) was not successful
before months.

\paragraph{Acknowledgements.}

Frédéric Chyzak and Manuel Kauers were supported in part
by the French-Austrian ANR-FWF grant EAGLES (ANR-22-CE91-0007 \& FWF-I-6130-N).
Frédéric Chyzak was further supported in part
by the European Research Council under the European Union's Horizon Europe research and innovation programme, grant agreement 101040794 (10000 DIGITS).
Hui Huang was supported in part by the NSFC grant (No.~12101105)
and the Natural Science Foundation of Fujian Province of China (No.~2024J01271).
The three authors were also supported in part
by the International Partnership Program of Chinese Academy of Sciences (Grant No.~167GJHZ2023001FN).
Part of the work was carried out when they were visiting the Academy of Mathematics and Systems Science in Beijing (China).

\printbibliography
\end{document}